\documentclass[preprints,article,accept,moreauthors,pdftex]{Definitions/mdpi} 

\setitemize{parsep=6pt,itemsep=0pt,leftmargin=*,labelsep=5.5mm,align=parleft}
\setenumerate{parsep=6pt,itemsep=0pt,leftmargin=*,labelsep=5.5mm,align=parleft}
\setlist[description]{itemsep=0mm}
\usepackage{bm, amssymb, tensor, tikz, makecell}
\newcommand{\T}{^{\mbox{\tiny T}}}
\newcommand{\Ts}{^{\,\mbox{\tiny T}}}
\newcommand{\B}[1]{{\bm #1}}

\newcommand{\ce}{constrained expression}
\newcommand{\ces}{constrained expressions}
\newcommand{\dd}[1]{\text{ d} #1}
\newcommand{\andd}{\quad \text{and} \quad}

\renewcommand{\C}[1]{\tensor*[^{}_{}]{\mathfrak{C}}{^{#1}_{}}}
\newcommand{\p}[2]{\tensor*[^{#1}_{}]{#2}{^{}_{}}}
\newcommand{\pC}[2]{\tensor*[^{#1}_{}]{\mathfrak{C}}{^{#2}_{}}}
\newcommand{\R}{\mathbb{R}}

\pdfoutput=1

\firstpage{1} 
\pubvolume{xx}
\issuenum{1}
\articlenumber{5}
\pubyear{2020}
\copyrightyear{2020}
\history{}

\Title{The Multivariate Theory of Functional Connections: Theory, Proofs, and Application in Partial Differential~Equations}

\Author{Carl Leake*, Hunter Johnston, and Daniele Mortari }
\AuthorNames{Carl Leake, Hunter Johnston, and Daniele Mortari}

\address[1]{Aerospace Engineering, Texas A\&M University, College Station, TX 77843, USA; hunterjohnston@tamu.edu (H.J.); mortari@tamu.edu (D.M.)}

\corres{\hangafter=1 \hangindent=1.05em \hspace{-0.82em} Correspondence: leakec@tamu.edu}

\abstract{This article presents a reformulation of the Theory of Functional Connections: a general methodology for \emph{functional interpolation} that can embed a set of user-specified linear constraints. The reformulation presented in this paper exploits the underlying functional structure presented in the seminal paper on the Theory of Functional Connections to ease the derivation of these interpolating functionals---called constrained expressions---and provides rigorous terminology that lends itself to straightforward derivations of mathematical proofs regarding the properties of these constrained expressions. Furthermore, the extension of the technique to and proofs in \texorpdfstring{$n$}{n}-dimensions is immediate through a recursive application of the univariate formulation. In all, the results of this reformulation are compared to prior work to highlight the novelty and mathematical convenience of using this approach. Finally, the methodology presented in this paper is applied to two partial differential equations with different boundary conditions, and, when data is available, the results are compared to state-of-the-art methods.}


\begin{document}

\section{Introduction}

The Theory of Functional Connections (TFC) is a mathematical framework used to construct functionals, functions of functions, that represent the family of all possible functions that satisfy some user-defined constraints; these functionals are referred to as ``\ces'' in the context of the TFC. In other words, the TFC is a framework for performing functional interpolation. In the seminal paper on TFC \cite{U-TFC}, a univariate framework was presented that could construct constrained expressions for constraints on the values of points or arbitrary order derivatives at points. Furthermore, Reference \cite{U-TFC} showed how to construct \ces\ for constraints consisting of linear combinations of values and derivatives at points, called linear constraints; for example, $y(x_1) + 3\pi \, y_{xx}(x_2) = 2 \, e$, for some points $x_1$ and $x_2$, where $y_{xx}$ symbolizes the second order derivative of $y$ with respect to $x$. In the current formulation, the univariate \ce\ has been used for a variety of applications, including solving linear and non-linear differential equations \cite{LDE, NDE}, hybrid systems \cite{hybrid_tfc}, optimal control problems \cite{EOL_EOI,FOL}, in quadratic and nonlinear programming \cite{QP_NLP}, and other applications \cite{Selected}. 

Recently, the TFC method has been extended to $n$-dimensions \cite{M-TFC}. This multivariate framework can provide functionals representing all possible $n$-dimensional manifolds subject to constraints on the value and arbitrary order derivative of $n-1$ dimensional manifolds. However, Reference \cite{M-TFC} does not discuss how the multivariate framework can be used to construct constrained expressions for linear constraints. Regardless, these multivariate \ces\ have been used to embed constraints into machine learning frameworks \cite{SVM,DeepTfc,XTFC} for use in solving partial differential equations (PDEs). Moreover, it was shown that this framework may be combined with orthogonal basis functions to solve PDEs \cite{M-TFC-PDE}; this is essentially the $n$-dimensional equivalent of the ordinary differential equations (ODEs) solved using the univariate formulation \cite{LDE,NDE}.

The contributions of this article are threefold. First, this article examines the underlying structure of univariate \ces\ and provides an alternative method for deriving them. This structure is leveraged to derive mathematical proofs regarding the properties of univariate \ces. Second, using the aforementioned structure, this article extends the multivariate formulation presented in Reference \cite{M-TFC} to include linear constraints by introducing the recursive application of univariate \ces\ as a method for generating multivariate \ces. Further, mathematical proofs are provided that prove the resultant constrained expressions indeed represent all possible manifolds subject to the given constraints. Thirdly, this article presents how the multivariate \ces\ can be combined with a linear expansion of $n$-dimensional orthogonal basis functions to numerically estimate the solutions of PDEs. While Reference \cite{M-TFC-PDE} showed that solving PDEs with the multivariate TFC is possible, it merely gave a cursory overview, skipping some rather important details; this article fills in those gaps.

The remainder of this article is structured as follows. Section \ref{section:univariate} introduces the univariate \ce, examines its underlying structure, and provides an alternative method to derive univariate \ces. Then, in Section \ref{sec:UnivariateProofs}, this structure is leveraged to rigorously define the univariate TFC \ce\ and provide some related mathematical proofs. In Section \ref{sec:Multivariate}, this new structure and the mathematical proofs are extended to $n$-dimensions, and a compact tensor form of the multivariate \ce\ is provided. Section \ref{sec:PDEApplication} discusses how to combine the multivariate \ce\ with multivariate basis functions to estimate the solutions of PDEs. Then, in Section \ref{sec:Results}, this method is used to estimate the solution of two PDEs, and the results are compared with state-of-the-art methods when data is available. Finally, Section \ref{sec:Conclusions} summarizes the article and provides some potential future directions for follow-on research.

\section{Univariate TFC}\label{section:univariate}

Extending the multivariate TFC to include linear constraints requires recursive applications of the univariate TFC. Hence, it is paramount the reader understand univariate TFC before moving to the multivariate case. First, the general form of the univariate \ce\ will be presented, followed by a few examples. These examples serve to solidify the readers understanding of the univariate \ce, as well as highlight nuances of deriving \ces. In addition, this section includes mathematical proofs that univariate TFC constrained expressions indeed represent the family of all possible functions that satisfy the constraints. 

Given a set of $k$ constraints, the univariate \ce\ takes the following form,
\begin{equation}\label{eq:uniCE}
    y (x, g (x)) = g (x) + \sum_{j = 1}^k s_j (x) \, \eta_j,
\end{equation}
where $g(x)$ is a free function, $s_j(x)$ are $k$ mutually linearly independent functions called support functions, and $\eta_j$ are $k$ coefficient functionals that are solved by imposing the constraints. The free function $g(x)$ can be chosen to be any function provided that it is defined at the constraints' locations. 

The following examples start from Equation~\eqref{eq:uniCE}, the framework proposed in the seminal paper on TFC \cite{U-TFC}, and highlight a unified structure that underlies the univariate TFC \ces. Following these examples is a section that rigorously defines this structure and provides important mathematical proofs.

\subsection{Univariate Example \# 1: Constraints at a Point}
Constraints at a point consist of constraints on the value at a point and constraints on a derivative at a point. Take for example the follow constraints,
\begin{equation*}
    y(0) = 1, \quad y_x(1) = 2, \quad y(2) = 3.
\end{equation*}

For this example, the support functions are chosen to be $s_1 = 1$, $s_2 = x^2$, and $s_3 = x^3$. Following Equation~\eqref{eq:uniCE} and imposing the three constraints leads to the simultaneous set of equations
\begin{align*}
    y(0) &= 1 = g(0) + \eta_1\\
    y_x(1) &= 2 = g_x(1) + 2\eta_2 + 3\eta_3\\
    y(2) &= 3 = g(2) + \eta_1 + 4 \eta_2 + 8 \eta_3.
\end{align*}

\noindent Solving this set of equations for the unknowns $\eta_j$ leads to the solution,
\begin{align*}
    \eta_1 &= 1 - g(0) \\
    \eta_2 &= \frac{10 - 3 g (0) + 3 g (2) - 8 g_x (1)}{4} \\
    \eta_3 &= \frac{g (0) - g (2) + 2 g_x (1)}{2}.
\end{align*}

\noindent Substituting the coefficient functionals back into Equation~\eqref{eq:uniCE} and simplifying yields,
\begin{equation}\label{eq:uniEx1Soln}
\begin{aligned}
    y(x,g(x)) =\ &g(x) + \frac{-2 x^3 + 3 x^2 + 4}{4} \Big(1 - g (0)\Big) + \Big(-x^3 + 2 x^2\Big) \Big(2 - g_x (1)\Big)\\
    &+\frac{2 x^3 - 3 x^2}{4}\Big(3 - g (2)\Big).
\end{aligned}
\end{equation}

\noindent It is simple to verify that regardless of how $g (x)$ is chosen, provided $g (x)$ exists at the constraint points, Equation~\eqref{eq:uniEx1Soln} always satisfies the given constraints. 

The support functions in the previous example were selected as $s_1 = 1$, $s_2 = x^2$, and $s_3 = x^3$. However, these support functions could have been any mutually linearly independent set of functions that permits a solution for the coefficient functionals $\eta_j$: to clarify the latter of these requirements, consider the same constraints with support functions $s_1 = 1$, $s_2 = x$, and $s_3 = x^2$. Then, the set of equations with unknowns $\eta_j$ is,
\begin{equation*}
    \begin{bmatrix} 1 & 0 & 0 \\ 0 & 1 & 2 \\ 1 & 2 & 4\end{bmatrix} \begin{Bmatrix} \eta_1 \\ \eta_2 \\ \eta_3\end{Bmatrix} = \begin{Bmatrix} 1 - g (0) \\ 2 - g_x (1) \\ 3 - g (2)\end{Bmatrix}.
\end{equation*}

Notice that when using these support functions the matrix that multiplies the coefficient functionals is singular. Thus, no solution exists, and therefore, the support functions $s_1 = 1$, $s_2 = x$, and $s_3 = x^2$ are an invalid set for these constraints.

Note that the matrix singularity does not depend on the free function. This means that the singularity arises when a linear combination of the selected support functions cannot be used to interpolate the constraints. Therefore, the singularity of the support function matrix is dependent on both the support functions chosen and the specific constraints to be embedded. This raises another important restriction on the expression of the support functions, not only must they be linearly independent, but they must constitute an interpolation model that is consistent for the specified~constraints.

Notice that each term, except the term containing only the free function, in the constrained expression is associated with a specific constraint and has a particular structure. To illustrate, examine the first constraint term from Equation~\eqref{eq:uniEx1Soln}, 
\begin{equation*}
    \underbrace{\frac{-2 x^3 + 3 x^2 + 4}{4}}_{\phi_1 (x)}\underbrace{(1 - g (0))}_{\rho_1(x,g(x))}.
\end{equation*}

The first term in the product, $\phi_1 (x)$, is called a {\it switching function}---Reference \cite{U-TFC} introduced these switching functions as ``coefficient'' functions, $\B{\beta}_k$---and is a function that is equal to $1$ when evaluated at the constraint it is referencing, and equal to $0$ when evaluated at all the other constraints. The second term of the product, $\rho_1 (x, g (x))$, is called a {\it projection functional}, and is derived by setting the constraint function equal to zero and replacing $y(x)$ with $g(x)$. In the case of constraints at a point, this is simply the difference between the constraint value and the free function evaluated at that constraint point. It is called the projection functional because it projects the free function to the set of functions that vanish at the constraint. When evaluating the switching function, $\phi_1 (x)$, at the constraint it is referencing it is equal to 1 (i.e., $\phi_1 (0) = 1$), and when it is evaluated at the other constraints it is equal to $0$ (i.e., $\frac{\partial \phi_1}{\partial x} (1) = 0$ and $\phi_1 (2) = 0$). The projection functional, $\rho_1 (x, g (x))$, is just the difference between the constraint $y (0) = 1$ and the free function evaluated at the constraint point, $g (0)$. This structure is important, as it shows up in the other constraint types too. Property \ref{prop:proj1} follows from the definition of the projection functional.

\begin{Property}\label{prop:proj1}
The projection functionals for constraints at a point are always equal to zero if the free function, $g(x)$, is selected such that it satisfies the associated constraint. 
\end{Property}

For example, if $g(x)$ is selected such that $g(0)=1$, then the first projection functional in this example becomes $\rho_1(x,g(x)) = 1-g(0) = 0$. Based on this structure, an alternate way to define the constrained expression, shown in Equation~\eqref{eq:uniCeAlt},  can be derived,
\begin{equation}\label{eq:uniCeAlt}
    y (x, g (x)) = g (x) + \sum_{j = 1}^k \phi_j (x) \, \rho_j(x,g(x)).
\end{equation}

The projection functionals are simple to derive, but the switching functions require some attention. From their definition, these functions must go to $1$ at their associated constraint and $0$ at all other constraints. Based on this definition, the following algorithm for deriving the switching functions is~proposed:
\begin{enumerate}
    \item Choose $k$ support functions, $s_k(x)$.
    \item Write each switching function as a linear combination of the support functions with unknown~coefficients.
    \item Based on the switching function definition, write a system of equations to solve for the unknown~coefficients. 
\end{enumerate}

To validate that this algorithm works, we will use the same constraints and support functions and rederive the \ce\ shown in Equation~\eqref{eq:uniEx1Soln}. Hence, $\phi_1 (x) = s_i (x) \, \alpha_{i1}$, $\phi_2 (x) = s_i (x) \, \alpha_{i2}$, and $\phi_3 (x) = s_i (x) \, \alpha_{i3}$, for some as yet unknown coefficients $\alpha_{ij}$. Note that in the previous mathematical expressions and throughout the remainder of the paper, the Einstein summation convention is used to improve readability. Now, the definition of the switching function is used to come up with a set of equations. For example, the first switching function has the three equations,
\begin{equation*}
    \phi_1(0) = 1, \quad \frac{\partial \phi_1}{\partial x}(1) = 0, \quad \text{and} \quad \phi_1(2) = 0.
\end{equation*}

\noindent These equations are expanded in terms of the support functions,
\begin{align*}
    \phi_1(0) &= (1) \cdot \alpha_{11} + (0) \cdot\alpha_{21} + (0) \cdot\alpha_{31} = 1\\
    \frac{\partial \phi_1}{\partial x}(1) &= (0)\cdot \alpha_{11} + (2) \cdot\alpha_{21} + (3) \cdot\alpha_{31} = 0\\
    \phi_1(2) &= (1) \cdot\alpha_{11} + (4) \cdot\alpha_{21} + (8) \cdot\alpha_{31} = 0,
\end{align*}
which can be compactly written as,
\begin{equation*}
    \begin{bmatrix} 1 & 0 & 0 \\ 0 & 2 & 3 \\ 1 & 4 & 8\end{bmatrix} \begin{Bmatrix} \alpha_{11} \\ \alpha_{21} \\ \alpha_{31} \end{Bmatrix} = \begin{Bmatrix} 1 \\ 0 \\ 0 \end{Bmatrix}.
\end{equation*}

\noindent The same is done for the other two switching functions to produce a set of equations that can be solved by matrix inversion.
\begin{align*}
    \begin{bmatrix} 1 & 0 & 0 \\ 0 & 2 & 3 \\ 1 & 4 & 8\end{bmatrix} \begin{bmatrix} \alpha_{11} & \alpha_{12} & \alpha_{13} \\ \alpha_{21} & \alpha_{22} & \alpha_{23} \\ \alpha_{31} & \alpha_{32} & \alpha_{33} \end{bmatrix} &= \begin{bmatrix} 1 & 0 & 0 \\ 0 & 1 & 0 \\ 0 & 0 & 1\end{bmatrix} \\
     \begin{bmatrix} \alpha_{11} & \alpha_{12} & \alpha_{13} \\ \alpha_{21} & \alpha_{22} & \alpha_{23} \\ \alpha_{31} & \alpha_{32} & \alpha_{33} \end{bmatrix} &= \begin{bmatrix} 1 & 0 & 0 \\ 0 & 2 & 3 \\ 1 & 4 & 8\end{bmatrix}^{-1} = \begin{bmatrix} 1 & 0 & 0 \\ \frac{3}{4} & 2 & -\frac{3}{4} \\ -\frac{1}{2} & -1 & \frac{1}{2}\end{bmatrix}
\end{align*}

\noindent Substituting the constants back into the switching functions and simplifying yields,
\begin{equation*}
    \phi_1 (x) = \frac{-2 x^3 + 3 x^2 + 4}{4}, \quad \phi_2 (x) =-x^3 + 2 x^2, \quad \text{and} \quad \phi_3 (x) = \frac{2 x^3 - 3 x^2}{4}.
\end{equation*}

\noindent Substituting the projection functionals and switching functions back into the constrained expression~yields,
\begin{equation*}
    y (x,g(x)) = g (x) + \frac{-2 x^3 + 3 x^2 + 4}{4}\Big(1 - g (0)\Big) + \Big(-x^3 + 2 x^2\Big)\Big(2 - g_x (1)\Big) + \frac{2 x^3 - 3 x^2}{4}\Big(3 - g (2)\Big),
\end{equation*}
which is identical to Equation~\eqref{eq:uniEx1Soln}. This approach to derive \ces\ using switching functions has, similar to the first approach, the risk of obtaining a singular matrix if the support functions selected are not able to interpolate the constraints. However, as will be demonstrated in the coming sections, this approach can be easily extended to multivariate domains via recursive applications of the univariate theory, and this approach lends itself nicely to mathematical proofs.

\subsection{Univariate Example \# 2: Linear Constraints}\label{sect:linear}

Linear constraints consist of linear combinations of the previous types of constraints. Note that by this definition, relative constraints such as $y(0) = y(1)$ are just a special case of linear constraints. Take for example the following two constraints,
\begin{equation*}
    y (0) = y (1), \qquad \text{and} \qquad 3 = 2 y (2) + \pi y_{xx} (0).
\end{equation*}

To generate a constrained expression, the projection functionals and switching functions must be found. Similar to the constraints at a point, first the constraints are arranged such that one side of the constraint is equal to zero; for example,
\begin{equation*}
    y (0) - y (1) = 0 \qquad \text{and} \qquad 3 - 2 y (2) - \pi y_{xx} (0) = 0.
\end{equation*}

\noindent Then, the projection functionals can be defined by replacing $y (x)$ with $g (x)$. Thus,
\begin{equation*}
    \rho_1 (x, g(x)) = g (0) - g (1) \andd \rho_2(x, g(x)) = 3 - 2 g (2) - \pi g_{xx} (0).
\end{equation*}

The switching functions are again defined such that they are equal to $1$ when evaluated with their associated constraint, and equal to $0$ when evaluated at all other constraints. The word ``evaluation'' in the previous sentence requires clarification. Substitution of the constrained expression back into the constraint should result in the expression $0 = 0$ (i.e., the constraint is satisfied). When doing so, the switching functions, $\phi (x)$, will be evaluated in the same way $y (x)$ is evaluated in the constraint. Thus, the constants within the constraint are not used in the evaluation. Moreover, because the projection functional is designed to exactly cancel the values of the free function in the constraint, the switching function equations should have the opposite sign. Hence, evaluation means to replace the function with the switching function, remove any terms not multiplied by the switching function, and multiply the entire equation by $-1$. Any reader confused by this linguistic definition of switching function evaluation may refer to Property \ref{prop:SwitchDefinition}, which defines the switching function evaluation mathematically. For this example, this leads to,
\begin{equation*}
    \phi_1 (1) - \phi_1 (0) = 1, \qquad 2 \phi_1 (2) + \pi \frac{\partial^2 \phi_1}{\partial x^2}(0) = 0,
\end{equation*}
for the first switching function, and 
\begin{equation*}
     \phi_2(1)-\phi_2(0) = 0, \qquad 2 \phi_2(2)+\pi \frac{\partial^2 \phi_2}{\partial x^2}(0) = 1,
\end{equation*}
for the second switching function. Note that while this ``evaluation'' definition may seem convoluted at first, it is in fact exactly what was done for the constraints at a point case. However, in that case, due to the simple nature of the constraints and the way the projection functionals were defined, this was simply the switching function evaluated at the point. 

Similar to the constraints at a point case, the switching functions are defined as a linear combination of support functions with unknown coefficients. Again, this can be written compactly in matrix form. For this example, the support functions $s_1(x) = 1$ and $s_2(x) = x$ are chosen. Then,
\begin{align*}
    \begin{bmatrix} 0 & 1 \\ 2 & 4 \end{bmatrix} \begin{bmatrix} \alpha_{11} & \alpha_{12}\\  \alpha_{21} & \alpha_{22} \end{bmatrix} &= \begin{bmatrix} 1 & 0 \\ 0 & 1 \end{bmatrix}\\
    \begin{bmatrix} \alpha_{11} & \alpha_{12}\\  \alpha_{21} & \alpha_{22} \end{bmatrix} &= \begin{bmatrix} 0 & 1 \\ 2 & 4 \end{bmatrix}^{-1} = \begin{bmatrix} -2 & \frac{1}{2} \\ 1 & 0 \end{bmatrix},
\end{align*}
which results in the switching functions,
\begin{equation*}
    \phi_1(x) = x-2, \quad \phi_2(x) = \frac{1}{2}.
\end{equation*}

Substituting the switching and projection functionals back into the constrained expression form given in Equation~\eqref{eq:uniCeAlt} yields,
\begin{equation*}
    y(x,g(x)) = g(x) + (x-2)\Big(g(0)-g(1)\Big)+\frac{1}{2}\Big(3-2g(2)-\pi g_{xx}(0)\Big).
\end{equation*}

\noindent By substituting this expression for $y$ back into the constraints, one can verify that indeed this constraint expression satisfies the constraints regardless of the choice of free function $g(x)$. Property \ref{prop:proj2} extends Property \ref{prop:proj1} to linear constraints.

\begin{Property}\label{prop:proj2}
The projection functionals for linear constraints are always equal to zero if the free function is selected such that it satisfies the associated constraint.
\end{Property}
\noindent For example, if $g(x)$ is selected such that $g(1) = g(0)$, then the first projection functional in this example becomes $\rho_1(x,g(x)) = g(1)-g(0) = 0$.

\section{General Formulation of the Univariate TFC}\label{sec:UnivariateProofs}
This section rigorously defines the TFC \ce\ and provides some relevant proofs. First, a functional is defined and its properties are investigated.
\begin{Definition}
A functional, e.g., $f(x,g(x))$, has independent variable(s) and function(s) as inputs, and produces a function as an output.
\end{Definition}

\noindent Note that a function as defined here is coincident with the computer science definition of a functional. One can think of a functional as a map for functions. That is, the functional takes a function, $g(x)$, as its input and produces a function, $h(x) = f(x,g(x))$ for any specified $g(x)$, as its output. Since this article is focused on constraint embedding, or in other words functional interpolation, we will not concern ourselves with the domain/range of the input and output functions. Rather, we will discuss functionals only in the context of their potential input functions, hereon referred to as the domain of the functional, and potential output functions, hereon referred to as the codomain of the functional. 

\begin{figure}[H]
    \centering\includegraphics[width=.65\linewidth]{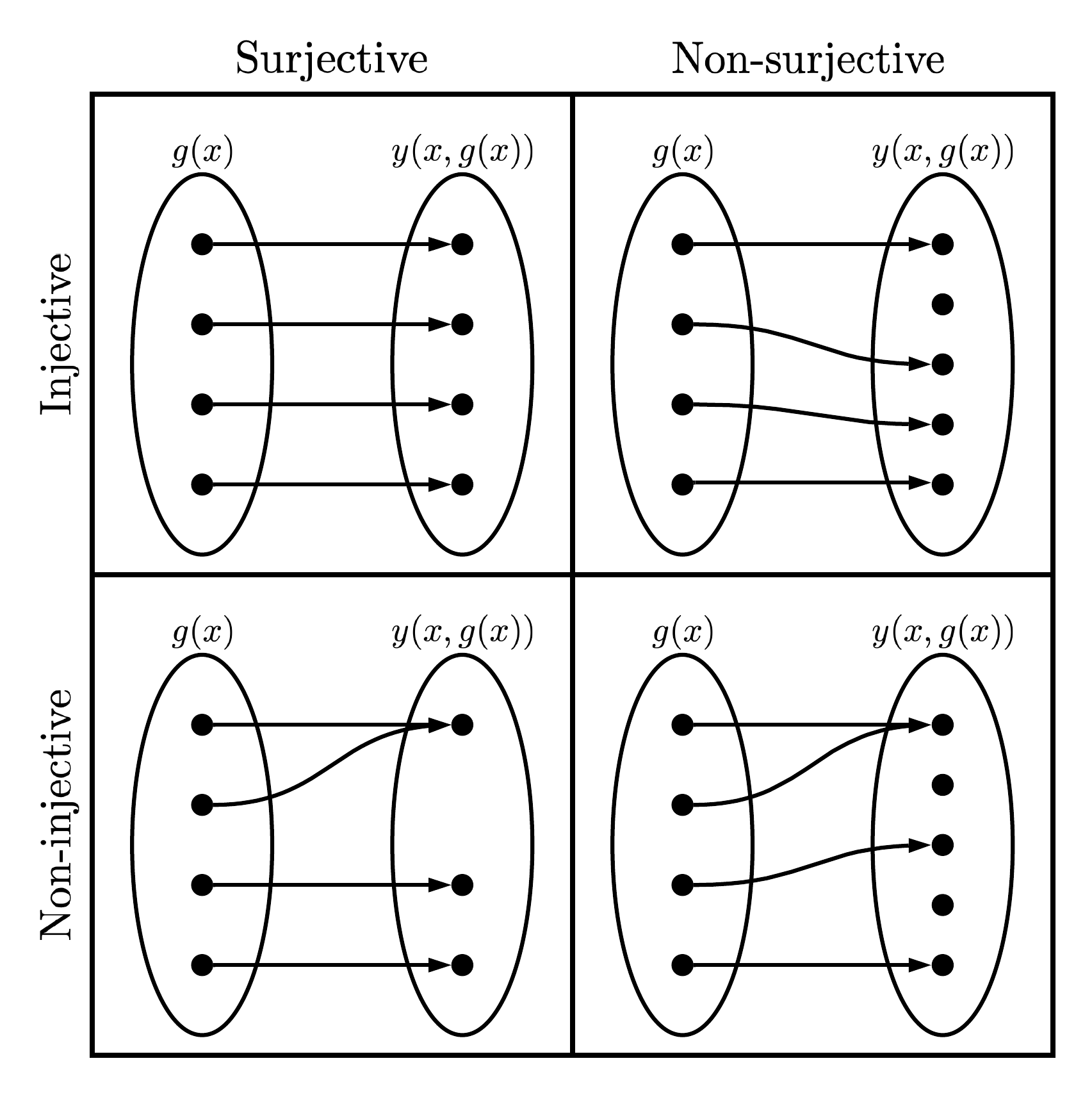}
    \caption{Graphical representation of injective and surjective functionals.}
    \label{fig:Injective_Surjective}
\end{figure}

Next, the definitions of injective, surjective, and bijective are extended from functions to~functionals. 
\begin{Definition}
A functional, $f(x,g(x))$, is said to be injective if every function in its codomain is the image of at most one function in its domain. 
\end{Definition}
\begin{Definition}
A functional, $f(x,g(x))$, is said to be surjective if for every function in the codomain, $h(x)$, there exists at least one $g(x)$ such that $h(x) = f(x,g(x))$.
\end{Definition}
\begin{Definition}
A functional, $f(x,g(x))$, is said to be bijective if it is both injective and surjective. 
\end{Definition}
To elaborate, Figure \ref{fig:Injective_Surjective} gives a graphical representation of each of these functionals, and examples of each of these functionals follow. Note that the phrase ``smooth functions'' is used here to denote continuous, infinitely differentiable, real valued functions. Consider the functional $f(x,g(x)) = e^{-g(x)}$ whose domain is all smooth functions and whose codomain is all smooth functions. The functional is injective because for every $h(x)$ in the codomain there is at most one $g(x)$ that maps $f(x,g(x))$ to $h(x)$.

However, the functional is not surjective, because the functional does not span the space of the codomain. For example, consider the desired output function $h(x) = -2$: there is no $g(x)$ that produces this output. Next, consider the functional $f(x,g(x)) = g(x)-g(0)$ whose domain is all smooth functions and whose codomain is the set of all smooth functions $h(x)$ such that $h(0) = 0$. This functional is surjective because it spans the space of all smooth functions that are $0$ when $x=0$, but it is not injective. For example, the functions $g(x) = x$ and $g(x) = x+3$ produce the same result, i.e., $f(x,x) = f(x,x+3) = x$. Finally, consider the functional $f(x,g(x)) = g(x)$ whose domain is all smooth functions and whose codomain is all smooth functions. This functional is bijective, because it is both injective and surjective. 

In addition, the notion of projection is extended to functionals. Consider the typical definition of a projection matrix $P^n = P$ for some $n\in\mathbb{Z}^+$. In other words, when $P$ operates on itself, it produces itself: a projection property for functionals can be defined similarly. 

\begin{Definition}
A functional is said to be a projection functional if it produces itself when operating on itself.
\end{Definition}
For example, consider a functional operating on itself, $f(x,f(x,g(x)))$. Then, if $f(x,f(x,g(x))) = f(x,g(x))$, then the functional is a projection functional. Note that proving $f(x,f(x,g(x))) = f(x,g(x))$ automatically extends to a functional operating on itself $n$ times: for example, $f(x,f(x,f(x,g(x))) = f(x,f(x,g(x))) = f(x,g(x))$, and so on.

Now that a functional and some properties of a functional have been investigated, the notation used in the prior section can be leveraged to rigorously define TFC related concepts. First, it is useful to define the constraint operator, denoted by the symbol $\C{}$.
\begin{Definition}
The constraint operator, $\C{i}$, is a linear operator that when operating on a function returns the function evaluated at the $i$-th specified constraint.
\end{Definition}

\noindent As an example, consider the 2nd linear constraint ($i = 2$) given in Section \ref{sect:linear}, $3 = 2 y(2) + \pi y_{xx}(0)$. For this problem, it follows that,
\begin{equation*}
    \C{2} [y(x)] =  2 y(2) + \pi y_{xx}(0).
\end{equation*}

\noindent The constraint operator is a linear operator, as it satisfies the two properties of a linear operator: (1)~$\C{i} [f(x) + g(x)] = \C{i}[f(x)] + \C{i}[g(x)]$ and (2) $\C{i}[a g(x)] = a\C{i}[g(x)]$. For example, again consider the 2$^\text{nd}$ linear constraint given in Section \ref{sect:linear},
\begin{align*}
    \C{2} [f(x)+g(x)] &= \C{2} [f(x)] + \C{2}[g(x)] = 2 f(2) + \pi f_{xx}(0) + 2 g(2) + \pi g_{xx}(0) \\
    \C{2}[a f(x)] &= a \C{2} [f(x)] = a \Big( 2 f(2) + \pi f_{xx}(0)\Big).
\end{align*}

\noindent Naturally, the constraint operator has specific properties when operating on the support functions, switching functions, and projection functionals.

\begin{Property}\label{prop:co_on_s}
The constraint operator acting on the support functions $s_j (x)$ produces the matrix \begin{equation*}
    \mathbb{S}_{ij} = \C{i}[s_j(x)].
\end{equation*}
\end{Property}
Again, consider the example from Section \ref{sect:linear} where the support functions were $s_1(x) = 1$ and $s_2(x) = x$. By applying the constraint operator,
\begin{equation*}
    \mathbb{S}_{ij} = \C{i}[s_j(x)] = \begin{bmatrix}\C{1}[s_1(x)] & \C{1}[s_2(x)] \\ \C{2}[s_1(x)] & \C{2}[s_2(x)] \end{bmatrix} = \begin{bmatrix}s_1(1) - s_1(0) & s_2(1) - s_2(0) \\ 2s_1(2) + \pi s_{1_{xx}}(0) & 2s_2(2) + \pi s_{2_{xx}}(0) \end{bmatrix} = \begin{bmatrix}0 & 1 \\ 2 & 4 \end{bmatrix},
\end{equation*}
which is identical to the matrix derived in Section \ref{sect:linear}. In fact, the matrix $\mathbb{S}_{ij}$ is simply the matrix multiplying the $\alpha_{ij}$ matrix in all the previous examples. Therefore, it follows that, $\mathbb{S}_{ij} \, \alpha_{jk} =  \alpha_{ij} \, \mathbb{S}_{jk} = \delta_{ik}$, where $\delta_{ik}$ is the Kroneker delta, and the solution of the $\alpha_{ij}$ coefficients are precisely the inverse of the constraint operator operating on the support functions. 

\begin{Property}\label{prop:SwitchDefinition}
The constraint operator acting on the switching functions $\phi_j(x)$ produces the Kronecker delta.  
\begin{equation*}
    \C{i}[\phi_j(x)] = \delta_{ij}
\end{equation*}
\end{Property}

\noindent This property is just a mathematical restatement of the linguistic definition of the switching function given earlier. One can intuit this property from the switching function definition, since they evaluate to $1$ at their specified constraint condition (i.e., $i=j$) and to $0$ at all other constraint conditions (i.e., $i \neq j$).

Using this definition of the constraint operator, one can define the projection functional in a compact and precise manner.
\begin{Definition}\label{def:projection_function}
Let $g(x)$ be the free function where $g: \R \to \R$, and let $\kappa_i\in\mathbb{R}$ be the numerical portion of the $i^{th}$ constraint. Then,
\begin{equation*}
    \rho_i(x,g(x)) = \kappa_i - \C{i}[g(x)]
\end{equation*}
\end{Definition}

\noindent Following the example from Section \ref{sect:linear}, the projection functional for the second constraint is,
\begin{equation*}
    \rho_2(x,g(x)) = \kappa_2 - \C{2}[g(x)] = 3 - 2 g(2) - \pi g_{xx}(0).
\end{equation*}

\noindent Note that in the univariate case $\kappa_i$ is a scalar value, but in the multivariate case $\kappa_i$ is a function. 


\begin{Theorem}\label{thrm:UniCe}
For any function, $f (x)$, satisfying the constraints, there exists at least one free function, $g (x)$, such that the TFC \ce\ $y(x,g(x)) = f(x)$.
\end{Theorem}

\begin{proof}
As highlighted in Properties \ref{prop:proj1} and \ref{prop:proj2}, the projection functionals are equal to zero whenever $g(x)$ satisfies the constraints. Thus, if $g(x)$ is a function that satisfies the constraints, then the \ce\ becomes $y (x, g (x)) = g (x) + \rho_i (x, g(x)) \, \phi_i (x) = g(x) + 0_i \, \phi_i (x) = g (x)$. Hence, by choosing $g (x) = f (x)$, the \ce\ becomes $y (x, f (x)) = f (x)$. Therefore, for any function satisfying the constraints, $f(x)$, there exists at least one free function $g (x) = f (x)$, such that the constrained expression is equal to the function satisfying the constraints, i.e., $y (x, f (x)) = f (x)$. 
\end{proof}

\begin{Theorem}\label{thrm:ProjUni}
The TFC univariate \ce\ is a projection functional.
\end{Theorem}

\begin{proof}
To prove Theorem~\ref{thrm:ProjUni}, one must show that $y (x, y (x, g (x))) = y (x, g (x))$. By definition, the constrained expression returns a function that satisfies the constraints. In other words, for any $g (x)$, $y (x, g (x))$ is a function that satisfies the constraints. From Theorem~\ref{thrm:UniCe}, if the free function used in the \ce\ satisfies the constraints, then the \ce\ returns that free function exactly. Hence, if the \ce\ functional is given itself as the free function, it will simply return itself. 
\end{proof}

\begin{Theorem}\label{thrm:NonUniG}
For a given function, $f (x)$, satisfying the constraints, the free function, $g (x)$, in the TFC \ce\ $y(x,g(x)) = f(x)$ is not unique. In other words, the TFC \ce\ is a surjective functional. 
\end{Theorem}

\begin{proof}
Consider the free function choice $g (x) = f (x) + \beta_j \,  s_j (x)$ where $\beta_j$ are scalar values on $\R$ and $s_j (x)$ are the support functions used to construct the switching functions $\phi_i (x)$.
\begin{equation*}
    y (x) = g (x) + \phi_i (x) \, \rho_i (x, g (x)).
\end{equation*}

\noindent Substituting the chosen $g (x)$ yields,
\begin{equation*}
    y (x) = f (x) + \beta_j \, s_j (x) + \phi_i (x) \, \rho_i (x, f (x) + \beta_j \, s_j (x)).
\end{equation*}

\noindent Now, according to Definition \ref{def:projection_function} of the projection functional,
\begin{equation*}
    y (x) = f (x) + \beta_j \, s_j(x) + \phi_i (x) \Big(\kappa_i - \C{i} [f(x) + \beta_j \, s_j(x)]\Big).
\end{equation*}

\noindent Since the constraint operator $\C{i}$ is a linear operator,
\begin{equation*}
    y(x) = f(x) + \beta_j s_j(x) +  \phi_i(x)\Big(\kappa_i - \C{i}[f(x)] -  \C{i}[s_j(x)]\beta_j\Big).
\end{equation*}

\noindent Since $f (x)$ is defined as a function satisfying the constraints, then $\C{i} [f(x)] = \kappa_i$, and,
\begin{equation*}
    y(x) = f(x) + \beta_j s_j(x) -  \phi_i(x)\C{i}[s_j(x)]\beta_j.
\end{equation*}

\noindent Now, according to Property \ref{prop:co_on_s} of the constraint operator, and by decomposing the switching functions~$\phi_i$,
\begin{equation*}
    y(x) = f(x) + \beta_j \, s_j(x) -  \alpha_{ki} \, s_k(x) \mathbb{S}_{ij} \, \beta_j.
\end{equation*}

\noindent Collecting terms results in,
\begin{equation*}
    y(x) = f(x) + \beta_j \Big(\delta_{jk} - \alpha_{ki} \, \mathbb{S}_{ij}\Big) s_k(x).
\end{equation*}

\noindent However, $\mathbb{S}_{ki} \alpha_{ij} = \delta_{kj}$ because $\alpha_{ij}$ is the inverse of $\mathbb{S}_{ki}$. Therefore, by the definition of inverse,
$\mathbb{S}_{ki} \alpha_{ij} =  \alpha_{ki} \mathbb{S}_{ij}   = \delta_{kj}$, and thus,
\begin{equation*}
    y(x) = f(x) + \beta_j \Big(\delta_{jk} - \delta_{jk} \Big) s_k(x).
\end{equation*}

\noindent Simplifying yields the result,
\begin{equation*}
    y(x) = f(x),
\end{equation*}
which is independent of the $\beta_js_j(x)$ terms in the free function. Therefore, the free function is not~unique.
\end{proof}

Notice that the non-uniqueness of $g(x)$ depends on the support functions used in the \ce, which has an immediate consequence when using \ces\ in optimization. If any terms in $g (x)$ are linearly dependent to the support functions used to construct the constrained expression, their contribution is negated and thus arbitrary. For some optimization techniques it is critical that the linearly dependent terms that do not contribute to the final solution be removed, else, the optimization technique becomes impaired. For example, prior research focused on using this method to solve ODEs \cite{LDE,NDE} through a basis expansion of $g(x)$ and least-squares, and the basis terms linearly dependent to the support functions had to be omitted from $g(x)$ in order to maintain full rank matrices in the least-squares.

The previous proofs coupled with the functional and functional property definitions given earlier provide a more rigorous definition for the TFC \ce: the TFC \ce\ is a surjective, projection functional whose domain is the space of all real-valued functions that are defined at the constraints and whose codomain is the space of all real-valued functions that satisfy the constraints. It is surjective because it spans the space of all functions that satisfy the constraints, its codomain, based on Theorem~\ref{thrm:UniCe}, but is not injective, because Theorem~\ref{thrm:NonUniG} shows that functions in the codomain are the image of more than one function in the domain: the functional is thus not bijective either because it is not injective. Moreover, the TFC \ce\ is a projection functional as shown in Theorem~\ref{thrm:ProjUni}. 

\section{Multivariate TFC}\label{sec:Multivariate}

Consider the general multivariate function $F: \R^n \to \R^m$ where $F = (f_1, f_2, \cdots, f_m)$. In this definition, $F$ is composed of the real-values functions $f_i: \R^n \to \R$ such that $f_i = f_i(x_1, x_2, \cdots, x_n)$ where $x_i$ are the independent variables. In terms of the TFC, the functions $f_i$ can be expressed as individual constrained expressions, and therefore the extension to multidimensional functions only involves extending the original method developed in Section \ref{section:univariate} for $f(x)$ to $f(x_1, x_2, \cdots, x_n)$. Once completed, the extension to the original definition of $F$ is immediate: simply write a multivariate \ce\ for every $f_i$ in $F$.

In the following section, the multivariate TFC is developed using a recursive application of the univariate TFC. In this manner, it can be shown that this approach is a generalization of the original theory, and that all mathematical proofs for the univariate \ces\ can easily be extended to the multivariate \ces. Then, a tensor form of the multivariate \ce\ is introduced by simplifying the recursive method. The tensor formulation provides a succinct way to write multivariate \ces.

\subsection{Recursive Application of Univariate TFC}
As discussed above, our extension to the multivariate case is concerned with deriving the constrained expression for the form $f = f(x_1, x_2, \cdots, x_n)$. For a set of constraints in the multivariate case, one can first create the constrained expression all constraints on $x_1$ using the univariate TFC formulation. The resulting univariate constrained expression, which we can denote as, $\p{1}{f}$ is then used as the free function in a constrained expression that includes all the constraints on $x_2$ to produce the expression $\p{2}{f}$. This method carries on until the final independent variable, $x_n$, is reached and the expression $\p{n}{f} = f$ and is the multivariate constrained expression.

This concept is best shown through some simple examples. These examples have two spatial dimensions only one dependent variable (i.e., $F: \R^2 \to \R$) we adopt the following notation:
\begin{align*}
    F = f_1 &:= u \\
    (x_1,x_2) &:= (x, y).
\end{align*}
 
\subsubsection{Multivariate Example \# 1: Value and Derivative Constraints}
Take for example the following constraints in two dimensions,
\begin{equation*}
    u(0,y) = \sin(2 \pi y), \quad u_{x}(0,y) = 0, \quad u(x,0) = x^2, \andd u(x,1) = \cos(x)-1.
\end{equation*}

\noindent First, the constrained expression is built for the constraints involving $x$. Using, the univariate TFC, this can be written as,
\begin{equation*}
    \p{1}{u}(x,y,g(x,y)) = g(x,y) + \sin(2 \pi y) - g(0,y) - x g_{x}(0,y).
\end{equation*}

\noindent Then, $\p{1}{u}(x,y)$ is used as the free function in the constrained expression involving the constraints on $y$. Since this problem is two-dimensional, the resultant expression is the multivariate TFC constrained~expression.
\begin{equation}\label{eq:MulEx1Ce1}
    u(x,y,\p{1}{u}(x,y)) = \p{1}{u}(x,y) + (1-y)\Big(x^2-\p{1}{u}(x,0)\Big)+y\Big(\cos(x)-1-\p{1}{u}(x,1)\Big)
\end{equation}

\noindent Substituting $\p{1}{u}$ into Equation~\eqref{eq:MulEx1Ce1} and simplifying yields,
\begin{equation}\label{eq:MulEx1Ce}
\begin{aligned}
    u(x,y,g(x,y)) =\ &g(x,y)+\sin (2 \pi  y)-g(0,y)-x g_{x}(0,y)+(1-y) \Big(x^2-g(x,0)+g(0,0)\\&
    +xg_{x}(0,0)\Big)+y \Big(\cos (x)-1-g(x,1)+g(0,1)+ x g_{x}(0,1)\Big).
\end{aligned}
\end{equation}

Alternatively, one could first write the expression for the constraints on $y$,
\begin{equation*}
    \p{2}{u}(x,y,g(x,y)) = g(x,y) + (1-y)\Big(x^2-g(x,0)\Big)+y\Big(\cos(x)-1-g(x,1)\Big),
\end{equation*}
and use $\p{2}{u}(x,y)$ as the free function in a constrained expression for the constraints on $x$,
\begin{equation}\label{eq:MulEx1Ce2}
    u(x,y,\p{2}{u}(x,y)) = \p{2}{u}(x,y) + \sin(2 \pi y) - \p{2}{u}(0,y) - x \frac{\partial (\p{2}{u})}{\partial x}(0,y).
\end{equation}

\noindent Substituting $\p{2}{u}(x,y)$ into Equation~\eqref{eq:MulEx1Ce2} and simplifying yields,
\begin{align*}
    u(x,y,g(x,y)) =\ &g(x,y)+\sin (2 \pi  y)-g(0,y)-x g_{x}(0,y)+(1-y) \Big(x^2-g(x,0)+g(0,0)\\&
    +xg_{x}(0,0)\Big)+y \Big(\cos (x)-1-g(x,1)+g(0,1)+ x g_{x}(0,1)\Big),
\end{align*}
the exact same result as in Equation~\eqref{eq:MulEx1Ce}. Therefore, it does not matter in what order recursive univariate TFC is applied to produce multivariate constrained expressions, as the final result will be the same.

\begin{figure}[H]
    \centering\includegraphics[width=0.7\linewidth]{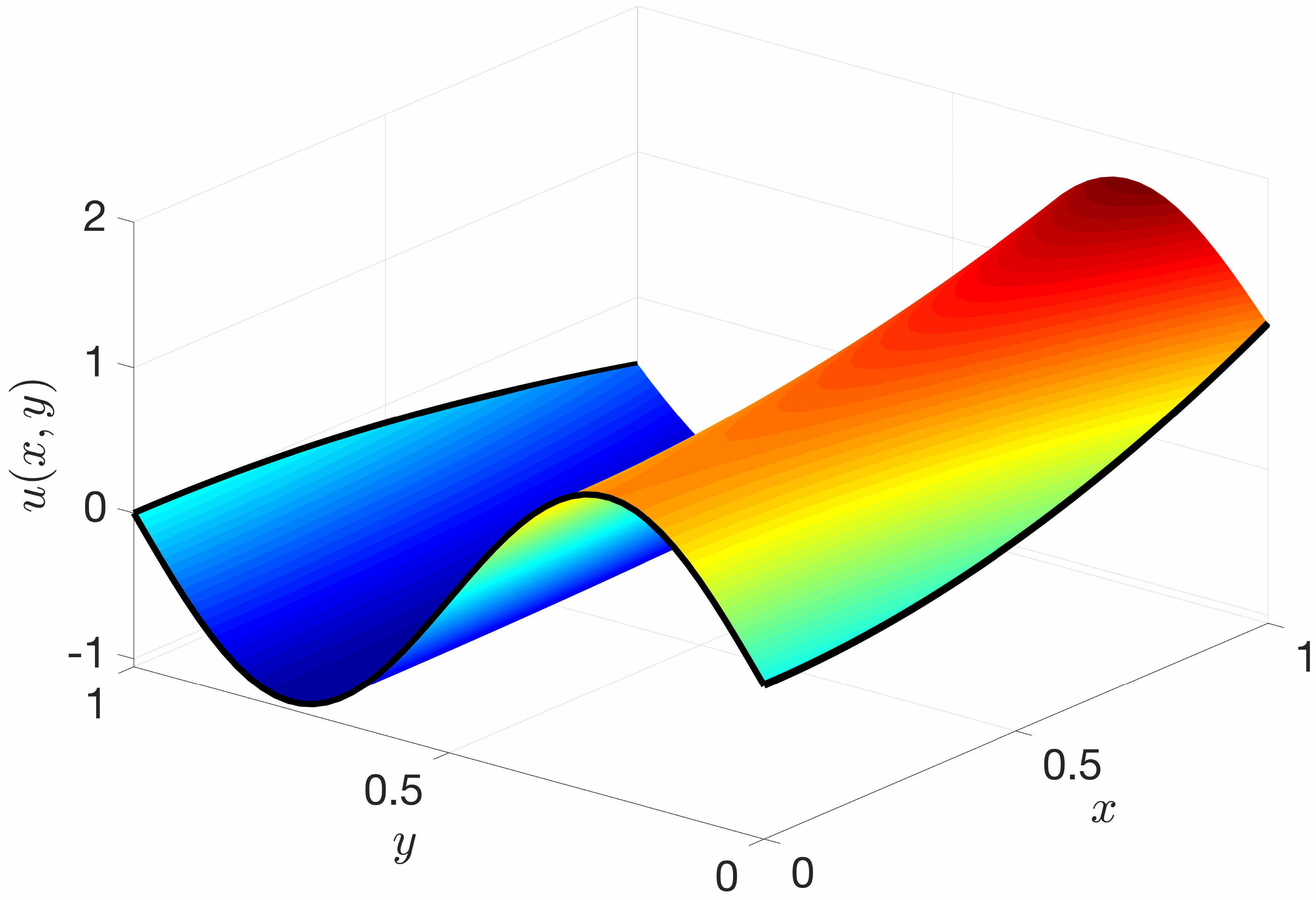}
    \caption{Multivariate example \# 1 constrained expression evaluated using $g(x,y) = x^2 cos(y) + 4$.}
    \label{fig:MultiEx1}
\end{figure}

Figure \ref{fig:MultiEx1} shows the constrained expression evaluated with the free function $g(x,y) = x^2 cos(y) + 4$. The constraints that can be visualized easily are shown as black lines. As expected, the TFC constrained expression satisfies these constraints exactly. 

\subsubsection{Multivariate Example \# 2: Linear Constraints}
Take for example the following constraints in two dimensions,
\begin{equation*}
   u(0,y) = y^2\sin(\pi y), \quad u(1,y)+u(2,y) = y\sin(\pi y), \quad u_y(x,0) = 0, \andd u(x,0) = u(x,1).
\end{equation*}

\noindent As in the first example, the univariate constrained expression is built for the constraints in $x$,
\begin{equation*}
    \p{1}{u}(x,y,g(x,y)) = g(x,y) + \frac{3-2x}{3}\Big(y^2\sin(\pi y)-g(0,y)\Big)+\frac{x}{3}\Big(\cos(\pi y)-g(2,y)-g(1,y)\Big).
\end{equation*}

\noindent Then, $\p{1}{u}(x,y)$ is used as the free function in the \ce\ for the constraints in $y$,
\begin{equation*}
    u(x,y,g(x,y)) = \p{1}{u}(x,y)-(y-y^2)\p{1}{u}_y(x,0)-y^2\Big(\p{1}{u}(x,1)-\p{1}{u}(x,0)\Big).
\end{equation*}

\noindent Substituting in $\p{1}{u}$ and simplifying yields,
\begin{equation}\label{eq:ceMultEx2}
\begin{aligned}
    u(x,y,\p{1}{u}(x,y)) =\ & g(x,y)+\left(y-y^2\right) \left(\frac{3-2 x}{3} g_y(0,0)-\frac{x}{3} \left(-g_y(1,0)-g_y(2,0)\right)-g_y(x,0)\right)\\
    &-y^2 \Big(\frac{3-2 x}{3} g(0,0)-\frac{3-2 x}{3} g(0,1)-\frac{x}{3}  (-g(1,0)-g(2,0))+\frac{x}{3}  (-g(1,1)\\
    &-g(2,1)) -g(x,0)+g(x,1)\Big)+\frac{3-2 x}{3} \left(y^2 \sin (\pi  y)-g(0,y)\right)\\
    &+\frac{x}{3}  \Big(-g(1,y) -g(2,y)+y \sin (\pi  y)\Big).
\end{aligned}
\end{equation}

\begin{figure}[H]
    \centering
    \includegraphics[width=0.6\linewidth]{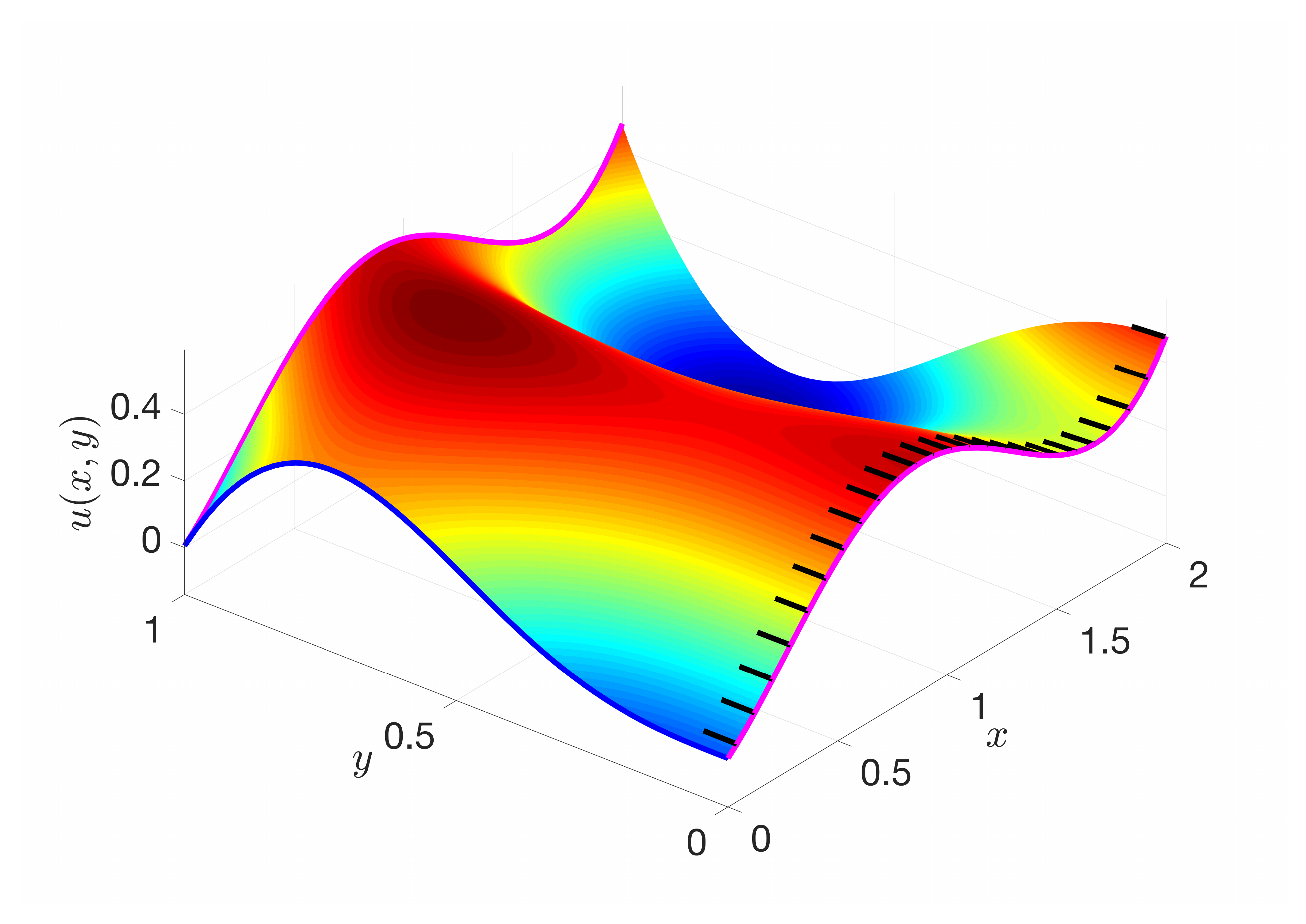}
    \caption{Multivariate example \# 2 constrained expression evaluated using $g (x, y) = x^2 \cos y + \sin (2 x)$. The blue line signifies the constraint on $u (0, y)$, the black lines signify the derivative constraint on $u_y (x, 0)$ and the magenta lines signify the relative constraint $u (x, 0) = u (x, 1)$. The linear constraint $u (1, y) + u (2, y) = y \sin(\pi y)$ is not easily visualized, but is nonetheless satisfied by the constrained~expression.}
    \label{fig:MultiEx2}
\end{figure}

Just as in the previous example, one could first write the constrained expression for the constraints in $y$, call it $\p{2}{u}(x,y)$, and then use $\p{2}{u}(x,y)$ as the free function in the constrained expression for the constraints in $x$: the result, after simplifying, would be identical to Equation~\eqref{eq:ceMultEx2}. Figure \ref{fig:MultiEx2} shows the constrained expression for the specific $g (x, y) = x^2 \cos y + \sin (2 x)$, where the blue line signifies the constraint on $u (0, y)$, the black lines signify the derivative constraint on $u_y (x, 0)$, and the magenta lines signify the relative constraint $u (x, 0) = u (x, 1)$. The linear constraint $u (1, y) + u (2, y) = y \sin(\pi y)$ is not easily visualized, but is nonetheless satisfied by the \ce.

\subsubsection{Multivariate Constrained Expression Theorems}
\begin{Theorem}\label{thrm:MultCe}
For any function, $f (\B{x})$, satisfying the constraints, there exists at least one free function, $g (\B{x})$, such that the multivariate TFC \ce\ $u(\B{x},g(\B{x})) = f(\B{x})$.
\end{Theorem}

\begin{proof}
Based on Theorem~\ref{thrm:UniCe}, the univariate constrained expression will return the free function if the free function satisfies the constraints. Let $ \p{1}{u}(\B{x})$ represent the univariate constrained expression for the independent variable $x_1$ that uses the free function $g(\B{x})$, $\p{2}{u}(\B{x})$ represent the univariate constrained expression for the independent variable $x_2$ that uses the free function $\p{1}{u}(\B{x})$, and so on up to $ \p{n}{u}(\B{x})$ which is simply the total constrained expression $u(\B{x})$. If we choose $g(\B{x}) = f(\B{x})$, then based on Theorem~\ref{thrm:UniCe} $\p{1}{u}(\B{x}) = f(\B{x})$. Applying Theorem~\ref{thrm:UniCe} recursively leads to $\p{2}{u}(\B{x}) = f(\B{x})$ and so on until $\p{n}{u}(\B{x}) = u(\B{x}) = f(\B{x})$. Hence, for any function satisfying the constraints, $f(\B{x})$, there exists a free function, $g(\B{x})=f(\B{x})$, such that the multivariate constrained expression is equal to the function satisfying the constraints, i.e., $u(\B{x},f(\B{x})) = f(\B{x})$.
\end{proof}

\begin{Theorem}\label{thrm:ProjMult}
The TFC multivariate \ce\ is a projection functional.
\end{Theorem}

\begin{proof}
To prove Theorem~\ref{thrm:ProjMult}, one must show that $u(\B{x},u(\B{x},g(\B{x}))) = u(\B{x},g(\B{x}))$. By definition, the constrained expression returns a function that satisfies the constraints. In other words, for any $g(\B{x})$, $u(\B{x},g(\B{x}))$ is a function that satisfies the constraints. From Theorem~\ref{thrm:MultCe}, if the free function used in the \ce\ satisfies the constraints, then the \ce\ returns that free function exactly. Hence, if the \ce\ function is given itself as the free function, it will simply return itself.
\end{proof}

\begin{Theorem}\label{thrm:NonMultG}
For a given function, $f (\B{x})$, satisfying the constraints, the free function, $g (\B{x})$, in the TFC \ce\ $u(\B{x},g(\B{x})) = f(\B{x})$ is not unique. In other words, the multivariate TFC \ce\ is a surjective functional. 
\end{Theorem}

\begin{proof}
Since each expression  $\p{i}{u}(\B{x})$ used in deriving the multivariate \ce\ is derived through the univariate formulation, then the results of the proof of Theorem~\ref{thrm:NonUniG} apply for each each $ \p{i}{u}(\B{x})$, and therefore the free function $g(\B{x})$ is not unique.
\end{proof}

Just like in the univariate case, this proof has immediate implications when using the constrained expression for optimization. Through the recursive application of the univariate TFC approach, any terms in $g(\B{x})$ that are linearly dependent to the the support functions, $s_i(x_1)$, $s_j(x_2)$, ... , $s_k(x_n)$, will not contribute to the solution. In the multivariate case, this also includes products of the support functions that include one and only one support function from each independent variable, e.g., $s_i(x_1)s_j(x_2)...s_k(x_n)$.

In addition, just as in the univariate case, Theorems \ref{thrm:MultCe}, \ref{thrm:ProjMult}, and \ref{thrm:NonMultG} allow for a more rigorous definition of the multivariate TFC \ce. The multivariate TFC \ce\ is a surjective, projection functional whose domain is the space of all real-valued functions that are defined at the constraints and whose codomain is the space of all real-valued functions that satisfy the~constraints.

\subsection{Tensor Form}
Recursive applications of the univariate TFC lead to expressions that lend themselves nicely to mathematical proofs, such as those in the previous section. However, for applications, it is typically more convenient to expression the constrained expression in a more compact form. Conveniently, the multivariate constrained expressions that are formed from recursive applications of the univariate TFC can be succinctly expressed in the following tensor form,
\begin{equation}
    u(\B{x}) = g(\B{x})+\mathcal{M}(\rho(\B{x},g(\B{x}))_{i_1i_2\dots i_n}\Phi_{i_1}(x_1)\Phi_{i_2}(x_2)\dots\Phi_{i_n}(x_n)
\end{equation}
where $i_1, i_2, \dots, i_n$ are $n$ indices associated with the $n$-dimensions that have constraints, $\mathcal{M}$ is an $n$-dimensional tensor whose elements are based on the projection functionals, $\rho(\B{x},g(\B{x}))$, and the $n$ vectors $\Phi$ are vectors whose elements are based on the switching functions for the associated~dimension. 

The $\mathcal{M}$ tensor can be constructed using a simple two-step process. Note that the arguments of functionals are dropped in this explanation for clarity. 
\begin{enumerate}
    \item The elements of the first order sub-tensors of $\mathcal{M}$ acquired by setting all but one index equal to one are a zero followed by the projection functionals for the dimension associated with that index. Mathematically,
    \begin{equation*}
        \mathcal{M}_{1\dots i_k \dots 1} = \begin{Bmatrix} 0, & \p{k}{\rho}_1, & \cdots, & \p{k}{\rho}_{\ell_k} \end{Bmatrix},
    \end{equation*}
    where $\p{k}{\rho}_j$ indicates the $j$-th projection functional of the $k$-independent variable and $\ell_k$ is the number of constraints associated with the $k$-th independent variable. 
    \item The remaining elements of the $\mathcal{M}$ tensor, those that have more than one index not equal to one, are the geometric intersection of the associated projection functionals multiplied by a sign ($-$~or~$+$). Mathematically, this can be written as,
    \begin{equation}\label{eq:GenMElement}
        \mathcal{M}_{i_1i_2 \dots i_n} = \pC{j}{i_j-1}\bigg[\pC{k}{i_k-1}\Big[ \cdots \big[\p{h}\rho_{i_h-1}\big] \cdots \Big]\bigg](-1)^m,
    \end{equation}
    where $i_j$, $i_k$, $\dots$, $i_h$ are the indices of $\mathcal{M}_{i_1i_2 \dots i_n}$ that are not equal to one and $m$ is equal to the number of non-one indices. If the constraint functions and free function are differential up to the order of derivatives required to compute Equation~\eqref{eq:GenMElement}, then by multiple applications of Clairaut's Theorem the constraint operators can be freely permuted \cite{M-TFC}. For example, Equation~\eqref{eq:GenMElement} could be re-written as,
    \begin{equation*}
        \mathcal{M}_{i_1i_2 \dots i_n} = \pC{h}{i_h-1}\bigg[\pC{j}{i_j-1}\Big[ \cdots \big[\p{k}\rho_{i_k-1}\big] \cdots \Big]\bigg](-1)^m.
    \end{equation*}
\end{enumerate}

The elements of the vectors $\Phi_{i_k}$ are simply composed of a $1$ followed by the switching functions associated with the $k$-th independent variable. Mathematically,
\begin{equation*}
    \Phi_{i_k} = \begin{Bmatrix} 1, & \p{k}{\phi}_1, & \cdots, & \p{k}{\phi}_{\ell_k}\end{Bmatrix},
\end{equation*}
where $\p{k}{\phi}_j$ denotes the $j$-th switching function of the $k$-th independent variable. 

To solidify the reader's understanding of the tensor form explained above, the constrained expressions for the two multivariate examples are re-derived using the tensor form. 

\subsubsection{Multivariate Example \# 1: Value and Derivative Constraints Using the Tensor Form}
The constraints from the first multivariate example are copied below for the reader's convenience. 
\begin{equation*}
    u(0,y) = \sin(2 \pi y), \quad u_{x}(0,y) = 0, \quad u(x,0) = x^2, \andd u(x,1) = \cos(x)-1
\end{equation*}

\noindent The projection functionals are defined based on the constraints,
\begin{equation*}
\begin{gathered}
    \p{1}{\rho}_1(x,y,g(x,y)) = \sin(2\pi y)-g(0,y),\quad \p{1}{\rho}_2(x,y,g(x,y)) = -g_x(0,y),\\
    \p{2}{\rho}_1(x,y,g(x,y)) = x^2-g(x,0), \andd \p{2}{\rho}_2(x,y,g(x,y)) = \cos(x)-1-g(x,1).
\end{gathered}
\end{equation*}

\noindent Then, the first step in constructing the $\mathcal{M}$ tensor can be completed.
\begin{equation*}
    \mathcal{M}_{ij}(x,y,g(x,y)) = \begin{bmatrix} 0 & x^2-g(x,0) & \cos(x)-1-g(x,1) \\ \sin(2\pi y)-g(0,y) & \text{-} & \text{-} \\ \cos(x)-1-g(x,1) & \text{-} & \text{-} \end{bmatrix}
\end{equation*}

\noindent The remaining elements of the $\mathcal{M}$ tensor are found in step 2 by calculating the geometric intersection of the projection functionals. For example, 
\begin{align*}
    \mathcal{M}_{22} &= \pC{1}{1}\Big[\p{2}{\rho}_1\Big](-1)^2 = -\Big[x^2-g(x,0)\Big]\Big|_{x=0} = g(0,0) \\ &=\pC{2}{1}\Big[\p{1}{\rho}_1\Big](-1)^2 = -\Big[\sin(2\pi y)-g(0,y)\Big]\Big|_{y=0} = g(0,0),
\end{align*}
where functional arguments have been dropped for clarity. The remaining elements are computed in a similar fashion to produce,
\begin{equation*}
    \mathcal{M}_{ij}(x,y,g(x,y)) = \begin{bmatrix} 0 & x^2-g(x,0) & \cos(x)-1-g(x,1) \\ \sin(2\pi y)-g(0,y) & g(0,0) & g(0,1) \\ \cos(x)-1-g(x,1) & g_x(0,0) & g_x(0,1) \end{bmatrix}.
\end{equation*}

\noindent The $\Phi$ vectors are assembled by concatenating a $1$ with the switching functions for that independent variable. Hence,
\begin{align*}
    \Phi_i(x) &= \begin{Bmatrix} 1 & 1 & x \end{Bmatrix} \\
    \Phi_j(y) &= \begin{Bmatrix} 1 & 1-y & y \end{Bmatrix}.
\end{align*}

Now, the tensor form of the constrained expression can be compactly written as,
\begin{equation*}
    u(x,y,g(x,y)) = g(x,y)+\mathcal{M}_{ij}(x,y,g(x,y)))\Phi_i(x)\Phi_j(y).
\end{equation*}

\noindent Note that expanding this expression produces the exact same constrained expression as the recursive~method.

\subsubsection{Multivariate Example \# 2: Linear Constraints Using the Tensor Form}
The constraints from the second multivariate example are copied below for the reader's~convenience. 
\begin{equation*}
   u(0,y) = y^2\sin(\pi y), \quad u(1,y)+u(2,y) = y\sin(\pi y), \quad u_y(x,0) = 0, \andd u(x,0) = u(x,1)
\end{equation*}

\noindent The projection functionals are defined based on the constraints,
\begin{equation*}
\begin{gathered}
    \p{1}{\rho}_1(x,y,g(x,y)) = y^2\sin(\pi y)-g(0,y), \quad \p{1}{\rho}_2(x,y,g(x,y)) = y\sin(\pi y)-g(1,y)-g(2,y),\\
    \p{2}{\rho}_1(x,y,g(x,y)) = -g_y(x,0), \andd \p{2}{\rho}_2(x,y,g(x,y)) = g(x,1)-g(x,0),
\end{gathered}
\end{equation*}
and the first step in constructing the $\mathcal{M}$ tensor is complete, 
\begin{equation*}
    \mathcal{M}_{ij}(x,y,g(x,y)) = \begin{bmatrix} 0 & -g_y(x,0) & g(x,1)-g(x,0) \\ y^2\sin(\pi y)-g(0,y) & \text{-} & \text{-} \\ y\sin(\pi y)-g(1,y)-g(2,y) & \text{-} & \text{-} \end{bmatrix}.
\end{equation*}

\noindent Then, just as in the previous example, remaining elements of the $\mathcal{M}$ tensor are found by calculating the geometric intersection of the projection functionals. For example,
\begin{align*}
    \mathcal{M}_{33} &= \pC{1}{2}\Big[\p{2}{\rho}_2\Big](-1)^2 = -\Big[g(x,1)-g(x,0)\Big]\Big|_{x=1}-\Big[g(x,1)-g(x,0)\Big]\Big|_{x=2} 
    \\&\quad\quad\quad\quad\quad\quad\quad\quad= g(1,0)-g(1,1)+g(2,0)-g(2,1) \\ &=\pC{2}{2}\Big[\p{1}{\rho}_2\Big](-1)^2 = \Big[y\sin(\pi y)-g(1,y)-g(2,y)\Big]\Big|_{y=1} -\Big[y\sin(\pi y)-g(1,y)-g(2,y)\Big]\Big|_{y=0} \\
    &\quad\quad\quad\quad\quad\quad\quad\quad= -g(1,1)-g(2,1)+g(1,0)+g(2,0),
\end{align*}
where functional arguments have been dropped for clarity. The complete $\mathcal{M}$ tensor for this example is,
\small 
\begin{equation*}
    \mathcal{M}_{ij}(x,y,g(x,y)) = \begin{bmatrix} 0 & -g_y(x,0) & g(x,1)-g(x,0) \\ y^2\sin(\pi y)-g(0,y) & g_y(0,0) & g(0,0)-g(0,1) \\ y\sin(\pi y)-g(1,y)-g(2,y) & g_y(1,0)+g_y(2,0) & g(1,0)-g(1,1)+g(2,0)-g(2,1) \end{bmatrix}.
\end{equation*}
\normalsize 

\noindent The $\Phi$ vectors are again assembled by concatenating a $1$ with the switching functions for the associated independent variable. 
\begin{align*}
    \Phi_i(x) &= \begin{Bmatrix} 1, & \frac{3-2x}{3}, & \frac{x}{3} \end{Bmatrix} \\
    \Phi_j(y) &= \begin{Bmatrix} 1, & y-y^2, & -y^2 \end{Bmatrix}
\end{align*}

Now, the tensor form of the constrained expression can be compactly written as,
\begin{equation*}
    u(x,y) = g(x,y)+\mathcal{M}_{ij}(x,y,g(x,y))\Phi_i(x)\Phi_j(y).
\end{equation*}

Note that expanding this expression produces the exact same constrained expression as the recursive~method.

\section{Applications to PDEs}\label{sec:PDEApplication}
In this article, orthogonal bases in $n$-dimensions, namely Chebyshev orthogonal polynomials of the first kind and Legendre orthogonal polynomials, are leveraged to approximate the solutions of PDEs with the TFC. For completeness, the equations to compute these polynomials are provided in Appendix \ref{appx:OP}. 

In general, multivariate basis sets can be created by using all possible products of the functions in the univariate basis sets. The measure that makes up the new multivariate basis set will be the product of measures of the univariate basis sets, and the domain of the multivariate basis set will be the union of the domains that make up the univariate basis sets. More details and insights regarding two-dimensional and $n$-dimensional orthogonal basis functions can be found in Refs. \cite{Ye,Dunkl,Xu}.

In this article, the free function will be defined as a linear combination of some multivariate basis set with unknown coefficients. The resultant constrained expression and its derivatives are substituted into the PDE. Since the free function consists of known basis functions and unknown coefficients, the PDE is transformed into an algebraic equation. This algebraic equation is discretized over the problem domain, and the unknown coefficients are used to minimize the residual of the PDE over the set of discretized points. The following subsections provide a detailed explanation of each major step, and a summary of the entire process is given in Figure \ref{fig:flowchart}.

\subsection{Defining the Free Function \texorpdfstring{$g(\B{x})$}{g(x)}}
Let us define $n$ independent variables in the vector $\B{x} = \{x_1, \cdots, x_k, \cdots, x_n\}\T$. Moreover, let the orthogonal basis set for each of these independent variables be denoted by $B_k^m$, where the superscript $m$ denotes the $m$-th basis function and the subscript $k$ denotes the $k$-th independent variable. For example, the third basis function for $x_2$ would be denoted as $B_2^3$. The domain of the multivariate basis will be denoted by $\Omega = \Omega_1 \times \Omega_2 \times \cdots \times \Omega_n$, where the generic $\Omega_k$ denotes the domain of the $k$-th basis set. Then, an arbitrary basis function for the multivariate domain can be written as,
\begin{equation}\label{eq:nDbasisAsTensorProduct}
    \mathcal{B} = B_1^{m_1} B_2^{m_2} \cdots B_n^{m_n},
\end{equation}
where $m_1,\cdots ,m_n \in \mathbb{Z^+}$. In other words, Equation~\eqref{eq:nDbasisAsTensorProduct} generates a multivariate basis via a tensor product of univariate basis functions \cite{nDbasisFunctions}. If one were to use all possible products of the functions in the individual basis sets, i.e., use all possible combinations of $m_1,\cdots ,m_n \in \mathbb{Z^+}$, an infinite set, then the resulting multivariate basis would span the union of the individual univariate basis sets' function spaces. However, when creating this expansion, one must pay attention to the results of Theorem~\ref{thrm:NonMultG}. Theorem~\ref{thrm:NonMultG} shows that the functions used to constrained the expression must be omitted from the formulation of $\mathcal{B}$. This ensures a full rank system in the later optimization steps.

As previously stated, the free function is chosen to be a linear combination of this multivariate basis with unknown coefficients. Mathematically, this can be expressed as,
\begin{equation}\label{eq:general_multi_g}
    g(\B{x}) =  \B{h}\Ts \B{\xi},
\end{equation}
where $\B{h}\in\mathbb{R}^{\left(\sum_{k=1}^n m_k\right)}$ whose elements are elements of $\mathcal{B}$, and $\B{\xi}$ is a same-sized vector of the unknown~coefficients.

\subsection{Derivatives of the Free Function}
In most applications, the domains $\Omega_k$ of the basis sets do not coincide with the domain of the problem (e.g., for Chebyshev and Legendre polynomials the functions are defined on $[-1, +1]$). Let the basis functions be defined on $z \in [z_0, z_f]$ and the problem be defined on $x_k \in [x_{k_0}, x_{k_f}]$ where $k$ corresponds to the dimension. In order to use the basis functions, a map between the basis function domain and problem domain must be created. The simplest map is a linear one,
\begin{equation}\label{eq:linearMapping}
z = z_0 + \frac{z_f-z_0}{x_{k_f}-x_{k_0}}(x - x_{k_0}) \quad \longleftrightarrow \quad x_k = x_{k_0} + \frac{x_{k_f}-x_{k_0}}{z_f-z_0}(z - z_0).
\end{equation}

\noindent The subsequent derivatives of the free function can be computed,
\begin{equation*}
    \frac{\partial^{n} g}{\partial x_k^{n}} = \left(\frac{\dd z}{\dd x_k}\right)^{n}  \frac{\partial^{n} \B{h}\Ts}{\partial z^n} \B{\xi},
\end{equation*}
but by defining,
\begin{equation}
c_k := \frac{\dd z}{\dd x_k} = \frac{z_f - z_0}{x_{k_f} - x_{k_0}},
\end{equation}
the derivative computations can be simply written as, 
\begin{equation*}
    \frac{\partial^{n} g}{\partial x_k^{n}} = c_k^{n} \frac{\partial^n \B{h}\Ts}{\partial z^n} \B{\xi}.
\end{equation*}

The immediate result is if the derivative of the function $g(\B{x})$ is taken with respect to the $x_k$ variable, then along with taking the derives of the basis functions with respect to this coordinate, the product must also be multiplied by the $c_k$ mapping coefficient. From this, it follows that a partial derivative with respect to multiple independent variables (e.g., $x_1$ and $x_2$) can be written as,
\vspace{-6pt}
\begin{equation*}
    \frac{\partial^{2} g}{\partial x_1 \partial x_2} = c_1 c_2 \frac{\partial \B{h}\Ts}{\partial x_1 \partial x_2} \B{\xi}.
\end{equation*}

\noindent This process applies to any derivative of the free function.

\subsection{Discretization}
In order to solve problems numerically, the problem domain (and therefore the basis function domain) must be discretized. Since this article uses Chebyshev and Legendre orthogonal polynomials, an optimal discretization scheme is the Chebyshev-Gauss-Lobatto nodes \cite{Collo1,Collo2}. For, $N+1$ points, the discrete points are calculated as,
\begin{equation}\label{eq:collo}
    z_j = -\cos\left(\frac{j \pi}{N}\right) \quad \text{for} \quad j = 0, 1, 2, \cdots, N.
\end{equation}

\noindent When compared with the uniform distribution, the collocation point distribution results in a much slower increase of the condition number of the matrix to be inverted in the least-squares as the number of basis functions, $m$, increases. The collocation points can be realized in the problem domain through the relationship provided in Equation~\eqref{eq:linearMapping}.

\subsection{Summary of the Major Steps to Solving PDEs}
\begin{figure}[H]
    \centering\includegraphics[width=\linewidth]{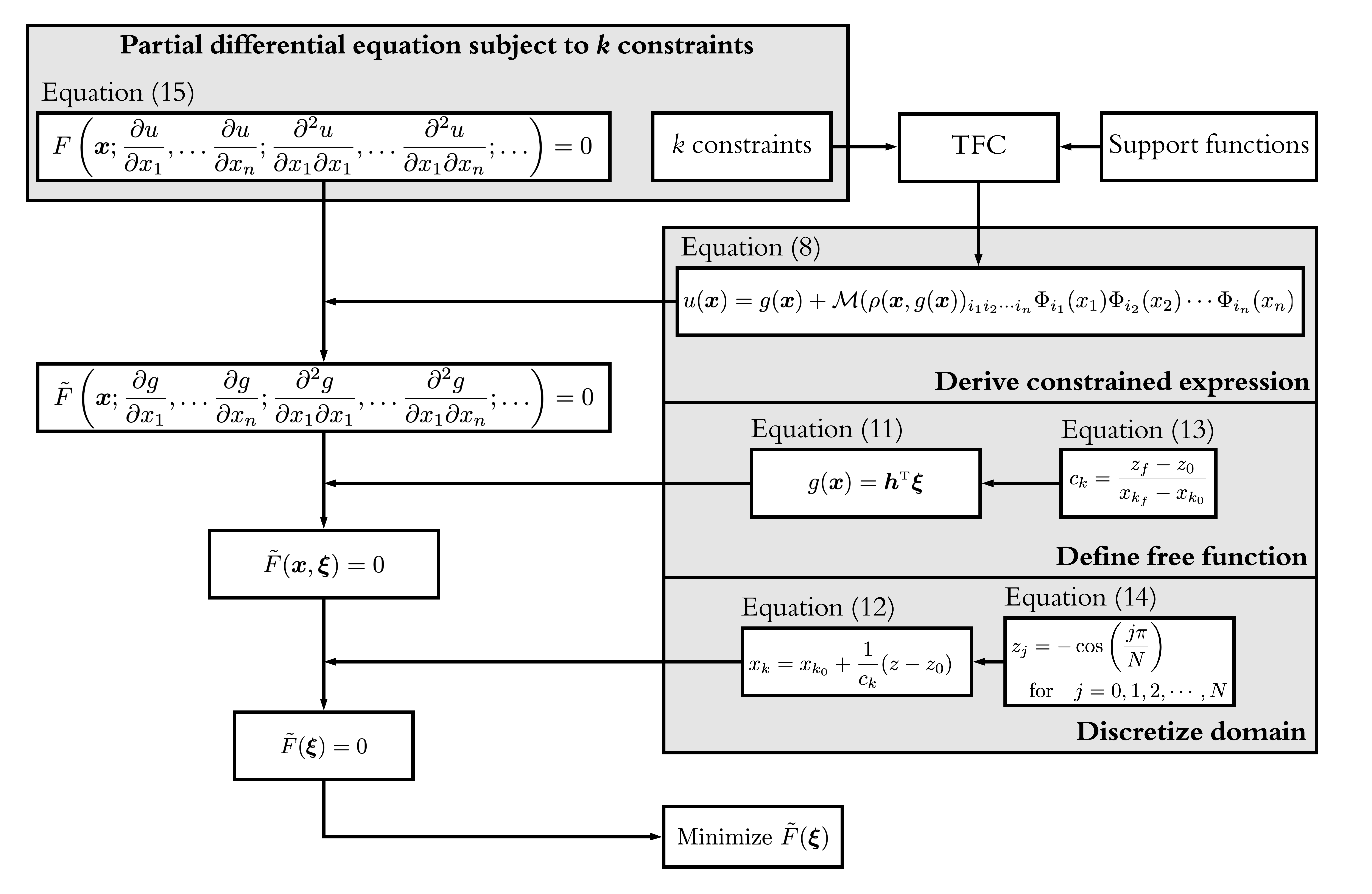}
    \caption{General flowchart of the TFC approach to solving partial differential equations.}
    \label{fig:flowchart}
\end{figure}

To summarize, consider a PDE for the function $u(\B{x})$ such that
\begin{equation}\label{eq:arbPDE}
    F\left(\B{x}; \frac{\partial u}{\partial x_{1}},\ldots {\frac {\partial u}{\partial x_{n}}};{\frac {\partial ^{2}u}{\partial x_{1}\partial x_{1}}},\ldots {\frac {\partial ^{2}u}{\partial x_{1}\partial x_{n}}};\ldots\right) = 0,
\end{equation}
subject to $k$ constraints. In general, the TFC approach to solving differential equations can be broken down into four major steps: (1) derive the constrained expression,  (2) define the free function,  (3) discretize the domain, and  (4) minimize the residual of the differential equation. The flow chart in Figure \ref{fig:flowchart} outlines these steps with all of the relevant equations.

In order to approximate the solution of Equation~\eqref{eq:arbPDE}, the \ce\ must be created: this is done using the theory developed in earlier sections. This constrained expression embeds the constraints of the differential equation, and by substituting the \ce\ into Equation~\eqref{eq:arbPDE}, transforms the original differential equation into an algebraic equation of the free function subject to no constraints. This transformed expression is denoted by the ($\sim$) symbol. Next, by defining the free function $g(\B{x})$ according to Equation~\eqref{eq:general_multi_g}, the differential equation becomes an algebraic equation in the unknown coefficient vector $\B{\xi}$. Lastly, this equation is discretized according to Equation~\eqref{eq:collo}, leading to a linear system of $\B{\xi}$ if the PDE is linear and a nonlinear system of $\B{\xi}$ if the PDE is nonlinear, which can be solved by many available optimization techniques. Previous work \cite{LDE,NDE} as well as the results provided in the following section utilize a least-squares approach.

\section{Results}\label{sec:Results}

This section applies the method to two PDEs. For each problem, the PDE and associated constraints are summarized along with the equations needed to construct the constrained expression. All numerical results were performed in Python, and utilized the autograd package \cite{autograd} to perform derivatives via automatic differentiation \cite{autodiff}. All computations were performed on a MacBook Pro (2016) macOS Version 10.15.3 with a 3.3 GHz Dual-Core Intel\textsuperscript{\textregistered} Core\texttrademark \, i7 and with 16 GB of RAM. All run times were calculated using the default\_timer function in the Python \verb"timeit" package. In all cases, matrix inversion was handled with NumPy's \verb"pinv" function.

Consequently, two specific computation times are provided (and tabulated in Appendix \ref{appx:data}), (1)~the full run-time of the problem and (2) the computational time associated with the least-squares. In these tables, the full run-time is drastically affected by the computational overhead from autograd, with full run times on the order of 0.5--50 seconds while the computation time for the least squares and nonlinear least squares is on the order of 0.5--185 milliseconds. 

For the results in the following sections, the accuracy of the method was determined according to the following process:
\begin{enumerate}
    \item Training
    \begin{enumerate}
        \item Estimate the solution of the PDE (i.e., determine the coefficients $\B{\xi}$) using the TFC method with $n$ points per independent variable and basis functions up to the to the $m$-th degree.
        \item Maximum training error: Using the training set discretization, determine the absolute error of the estimated solution compared with the true solution and record the maximum~value.
    \end{enumerate}
    \item Test
    \begin{enumerate}
        \item Using converged $\B{\xi}$ parameters from the training phase, refine the discretization of domain with $n = 100$ equally spaced points per dimension and evaluate estimated solution at these points.
        \item Maximum test error: Using the test set discretization, determine the absolute error of the estimated solution compared with the true solution and record the maximum value.
    \end{enumerate}
\end{enumerate}

Additionally, for both numerical tests, the method was completed over a varying range of discretization points per independent variable, $n$, and degree of basis expansion, $m$. 
The results in this section and in Appendix \ref{appx:data} are reported as a function of these two parameters. For example, a value of $n = 5$ would imply a $5\times 5$ grid or 25 points. Likewise, a value of $m = 5$ would imply that all the univariate basis functions, and combinations thereof, are at most quintic functions. However, the number of coefficients to be solved, and therefore the size of the matrix to be inverted, is dependent on the constrained expression, since some terms need to be removed from the expression of $g(x,y)$ (see Theorem~\ref{thrm:NonMultG}). The total number of basis functions associated with each degree for both Problems \#1 and \#2 are displayed in Table \ref{tab:numBase}.

\begin{table}[H]
\centering
\caption{Equivalence of number of basis function compared to degree of basis expansion for both Problems \#1 and \#2.}
\label{tab:numBase}
\begin{tabular}{cc} 
\toprule
 \makecell{\textbf{m}} & \makecell{\textbf{Number of} \\ \textbf{Functions}}\\  \midrule
5  & 17 \\
10 & 62 \\
15 & 132 \\
20 & 227 \\
25 & 347 \\
\bottomrule
\end{tabular}
\end{table}

\subsection{Problem \# 1}
Consider the PDE solved in Lagaris et al. \cite{OrigOdePde}, Mall \& Chakraverty \cite{CNN}, Sun et al. \cite{BNN}, and~Schiassi~et~al.~\cite{XTFC}, 
\begin{equation*}
    u_{xx}(x,y) + u_{yy} (x,y) = e^{-x}(x - 2 + y^3 + 6y)
\end{equation*}
where $x,y \in [0,1]$ and subject to,
\begin{eqnarray*}
    u(0,y) &=& y^3\\
    u(1,y) &=& (1+y^3)e^{-1}\\
    u(x,0) &=& xe^{-x}\\
    u(x,1) &=& e^{-x}(x+1),
\end{eqnarray*}
which has the true solution $u(x,y) = e^{-x}(x + y^3)$. Using the proposed method, the \ce\ can be derived and written in its tensor form as,
\begin{equation*}
    u(x,y,g(x,y)) = g(x,y)+\mathcal{M}_{ij}(x,y,g(x,y))\Phi_i(x)\Phi_j(y),
\end{equation*}
where $g(x,y)$ is defined according to Equation~\eqref{eq:general_multi_g} and for these numerical tests is implemented with either Chebyshev or Legendre polynomials. Furthermore, for this problem,
\begin{equation*}
    \mathcal{M}_{ij}(x,y,g(x,y)) = \begin{bmatrix} 0 & xe^{-x}-g(x,0) & e^{-x}(x+1)-g(x,1) \\
    y^3-g(0,y) & g(0,0) & -1+g(0,1) \\ 
    (1+y^3)e^{-1}-g(1,y) & -e^{-1}+g(1,0) & -2e^{-1}+g(1,1) \end{bmatrix}.
\end{equation*}
and
\begin{align*}
    \Phi_i(x) &= \begin{Bmatrix} 1, & 1-x, & x\end{Bmatrix}, \\
    \Phi_j(y) &= \begin{Bmatrix} 1, & 1-y, & y\end{Bmatrix}.
\end{align*}

\noindent It follows that the expanded \ce\ is,
\begin{align*}
    u(x,y) =\ &g(x,y)-(x-1) \left(y (-g(0,0)+g(0,1)-1)+g(0,0)+y^3\right)+(x-1) g(0,y)\\
    &+x (y g(1,1)-(y-1) g(1,0))-x g(1,y)+(y-1) g(x,0)-y g(x,1)+\frac{x y \left(y^2-1\right)}{e}\\
    &+e^{-x} (x+y).
\end{align*}

\begin{figure}[H]
    \centering
    \includegraphics[width=\linewidth]{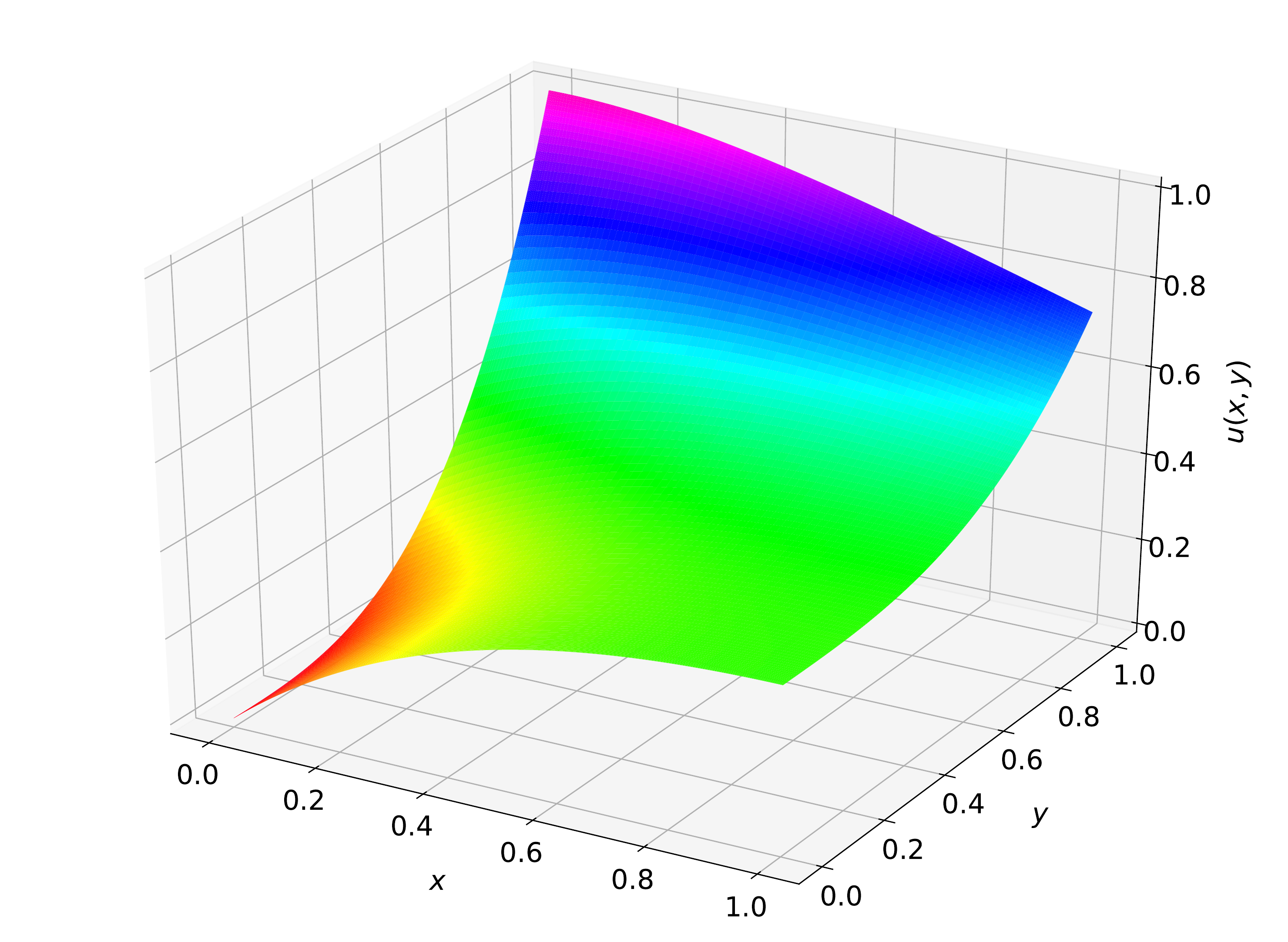}
    \caption{Problem \#1 analytical solution.}
    \label{fig:p1Soln}
\end{figure}

\noindent Figure \ref{fig:p1Soln} shows the analytical solution for this problem and Tables \ref{tab:p1_cheb} and \ref{tab:p1_leg} display the maximum error over the domain for the test set using both Chebyshev and Legendre polynomials.

It can be seen that the difference in accuracy of Chebyshev polynomials versus Legendre polynomials is negligible: the solutions that use Legendre polynomials maintain a slightly more accurate estimation than those that use Chebyshev polynomials when using higher order terms.

\begin{table}[H]
\caption{Maximum test set solution error using Chebyshev polynomials for Problem \# 1.}\label{tab:p1_cheb}
\[\begin{array}{c|ccccc}
\noalign{\hrule height 1.0pt}
\tikz{\node[below left, inner sep=1pt] (n) {\textbf{n}};%
      \node[above right,inner sep=1pt] (m) {\textbf{m}};%
      \draw (n.north west|-m.north west)--(n.south east-|m.south east);}

 & \textbf{5} & \textbf{10} & \textbf{15} & \textbf{20} & \textbf{25}\\
\hline
5  & 6.26\times 10^{-4} & \text{-}   & \text{-}   & \text{-}   & \text{-} \\
10 & 5.53\times 10^{-4} & 1.20\times 10^{-10} & \text{-}   & \text{-}   & \text{-} \\
15 & 5.30\times 10^{-4} & 1.17\times 10^{-10} & 4.44\times 10^{-16} & \text{-}   & \text{-} \\
20 & 5.20\times 10^{-4} & 1.16\times 10^{-10} & 5.00\times 10^{-16} & 4.44\times 10^{-16} & \text{-} \\
25 & 5.13\times 10^{-4} & 1.15\times 10^{-10} & 7.22\times 10^{-16} & 2.61\times 10^{-15} & 5.55\times 10^{-16}\\
30 & 5.09\times 10^{-4} & 1.14\times 10^{-10} & 6.66\times 10^{-16} & 8.88\times 10^{-16} & 3.22\times 10^{-15}\\  \noalign{\hrule height 1.0pt}
\end{array}\]
\end{table}
\unskip

\begin{table}[H]
\caption{Maximum test set solution error using Legendre polynomials for Problem \# 1.}\label{tab:p1_leg}
\[\begin{array}{c|ccccc}

\noalign{\hrule height 0.5pt}
\tikz{\node[below left, inner sep=1pt] (n) {\textbf{n}};%
      \node[above right,inner sep=1pt] (m) {\textbf{m}};%
      \draw (n.north west|-m.north west)--(n.south east-|m.south east);}
      
 & \textbf{5} & \textbf{10} & \textbf{15} & \textbf{20} & \textbf{25}\\
 
\hline
5  & 6.26\times 10^{-4} & \text{-} & \text{-} & \text{-} & \text{-} \\
10 & 5.53\times 10^{-4} & 1.20\times 10^{-10} & \text{-} & \text{-} & \text{-} \\
15 & 5.30\times 10^{-4} & 1.17\times 10^{-10} & 4.44\times 10^{-16} & \text{-}   & \text{-} \\
20 & 5.20\times 10^{-4} & 1.16\times 10^{-10} & 5.55\times 10^{-16} & 4.44\times 10^{-16} & \text{-} \\
25 & 5.13\times 10^{-4} & 1.15\times 10^{-10} & 4.44\times 10^{-16} & 4.44\times 10^{-16} & 5.55\times 10^{-16}\\
30 & 5.09\times 10^{-4} & 1.14\times 10^{-10} & 4.44\times 10^{-16} & 4.44\times 10^{-16} & 5.55\times 10^{-16}\\ \noalign{\hrule height 1.0pt}
\end{array}\]
\end{table}

Finally, the results of the numerical test for Problem \#1 are compared to the other approaches in Table \ref{tab:p1_comp}, where the maximum training and test errors are presented. It can be seen that this method produces an estimate at machine level precision that is at least 3 orders of magnitudes more accurate than the other methods. In fact, the next closest method is a TFC based approach where the free function is expressed using an Extreme Learning Machine \cite{ELM} in the paper by Schiassi et al. \cite{XTFC}.

\begin{table}[H]
\centering
\caption{Comparison of maximum training and test error of TFC with current state-of-the-art techniques for Problem \# 1.}
\label{tab:p1_comp}
\begin{tabular}{ccc} 
\toprule
\makecell{\textbf{Method}} & \makecell{\textbf{Training Set}\\\textbf{Maximum Error}} & \makecell{Test Set\\\textbf{Maximum Error}}\\   \midrule
{TFC} & {$2.22 \times 10^{-16}$} & {$4.44 \times 10^{-16}$} \\
{X-TFC \cite{XTFC}} & {$3.8 \times 10^{-13}$} & {$5.1 \times 10^{-13}$} \\
{FEM} & {$2 \times 10^{-8}$} & {$1.5 \times 10^{-5}$} \\
{Reference \cite{OrigOdePde}} & {$5 \times 10^{-7}$} & {$5 \times 10^{-7}$} \\
{Reference \cite{CNN}} & {\text{-}} & {$3.2 \times 10^{-2}$} \\
{Reference \cite{BNN}} & {\text{-}} & {$2.4 \times 10^{-4}$} \\
\bottomrule
\end{tabular}
\end{table}

\subsection{Problem \#2}
Problem \#2 is a PDE with a linear constraint created by the authors.
\begin{equation*}
    u_{xx}(x,y) + u_x (x,y) u_y(x,y) =2 \cos (y)-2 x^3 \sin (y) \cos (y)
\end{equation*}
where $(x,y) \in [0,1]\times[0,2\pi]$ and subject to,
\begin{eqnarray*}
    u(0,y) &=& 0\\
    u(1,y) &=& \cos(y)\\
    u(x,0) &=& u(x,1),
\end{eqnarray*}
which has the true solution $u(x,y) = x^2\cos(y)$.

Using the TFC, the \ce\ for these boundary conditions can be written in the tensor form as,
\begin{equation*}
    u(x,y,g(x,y)) = g(x,y)+\mathcal{M}_{ij}(x,y,g(x,y))\Phi_i(x)\Phi_j(y),
\end{equation*}
where
\begin{equation*}
    \mathcal{M}_{ij}(x,y,g(x,y)) = \begin{bmatrix} 0 & g(x,0)-g(x,2\pi)\\ 
    -g(0,y) & g(0,2\pi)-g(0,0)\\
    \cos(y)-g(1,y) & g(1,2\pi)-g(1,0)\end{bmatrix}
\end{equation*}
and
\begin{align*}
    \Phi_i(x) &= \begin{Bmatrix} 1 & 1-x & x \end{Bmatrix}\T,\\
    \Phi_j(y) &= \begin{Bmatrix} 1 & \frac{y}{2\pi} \end{Bmatrix}\T,
\end{align*}
or in its expanded form as,
\begin{align*}
    u(x,y,g(x,y)) =\ &g(x,y)-(1-x) \left(\frac{y (g(0,0)-g(0,2 \pi ))}{2 \pi }+g(0,y)\right)+\frac{y (g(x,0)-g(x,2 \pi ))}{2 \pi }\\
    &+x \left(-g(1,y)-\frac{y (g(1,0)-g(1,2 \pi ))}{2 \pi }+\cos (y)\right).
\end{align*}

\begin{figure}[H]
    \centering
    \includegraphics[width=\linewidth]{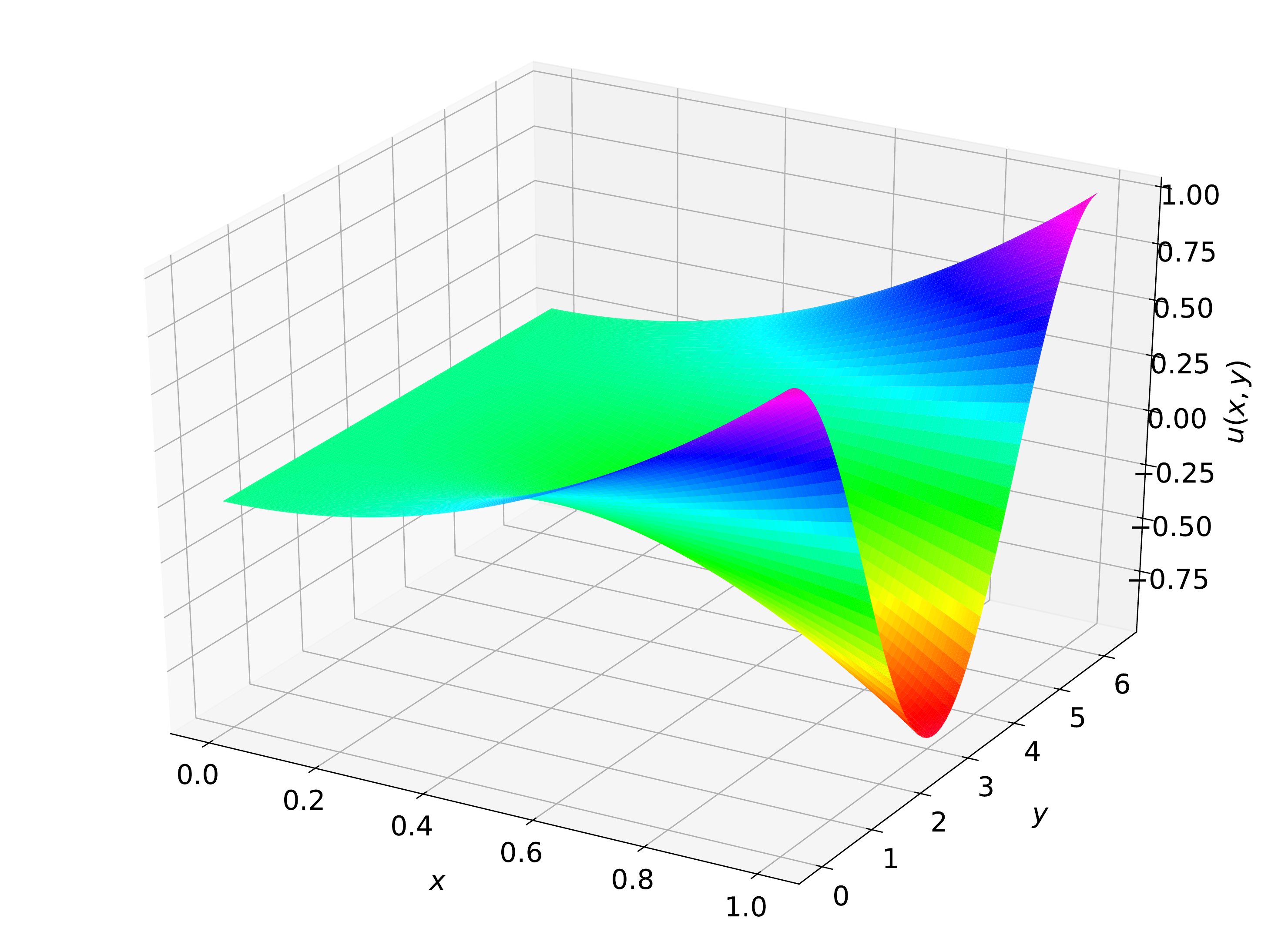}
    \caption{Problem \#2 analytical solution.}
    \label{fig:p2Soln}
\end{figure}

Figure \ref{fig:p2Soln} shows the analytical solution of Problem \#2, and Tables \ref{tab:p2_cheb} and \ref{tab:p2_leg} show the maximum test set error over the domain using Chebyshev and Legendre polynomials respectively. 

\unskip

\begin{table}[H]
\caption{Maximum test set solution error using Chebyshev polynomials for Problem \# 2.}\label{tab:p2_cheb}
\[\begin{array}{c|ccccc}
\noalign{\hrule height 1pt}

\tikz{\node[below left, inner sep=1pt] (n) {\textbf{n}};%
      \node[above right,inner sep=1pt] (m) {\textbf{m}};%
      \draw (n.north west|-m.north west)--(n.south east-|m.south east);}
 & \textbf{5} & \textbf{10} & \textbf{15} & \textbf{20} & \textbf{25}\\
\hline
5  & 1.03\times 10^{-1} & \text{-}   & \text{-}   & \text{-}   & \text{-} \\
10 & 9.20\times 10^{-2} & 2.49\times 10^{-5} & \text{-}   & \text{-}   & \text{-} \\
15 & 9.03\times 10^{-2} & 1.54\times 10^{-5} & 4.34\times 10^{-9} & \text{-}   & \text{-} \\
20 & 8.94\times 10^{-2} & 1.52\times 10^{-5} & 4.56\times 10^{-9} & 5.33\times 10^{-15} & \text{-} \\
25 & 8.88\times 10^{-2} & 1.50\times 10^{-5} & 4.53\times 10^{-9} & 2.72\times 10^{-15} & 4.44\times 10^{-16}\\
30 & 8.85\times 10^{-2} & 1.49\times 10^{-5} & 4.51\times 10^{-9} & 2.72\times 10^{-15} & 3.33\times 10^{-16}\\ \noalign{\hrule height 1pt}
\end{array}\]
\end{table}
\unskip

\begin{table}[H]
\caption{Maximum test set solution error using Legendre polynomials for Problem \# 2.}\label{tab:p2_leg}
\[\begin{array}{c|ccccc}

\noalign{\hrule height 1pt}

\tikz{\node[below left, inner sep=1pt] (n) {\textbf{n}};%
      \node[above right,inner sep=1pt] (m) {\textbf{m}};%
      \draw (n.north west|-m.north west)--(n.south east-|m.south east);}
 & \textbf{5} & \textbf{10} & \textbf{15} & \textbf{20} & \textbf{25}\\
\hline
5  & 1.03\times 10^{-1} & \text{-} & \text{-} & \text{-} & \text{-} \\
10 & 9.20\times 10^{-2} & 2.49\times 10^{-5} & \text{-} & \text{-} & \text{-} \\
15 & 9.03\times 10^{-2} & 1.54\times 10^{-5} & 4.34\times 10^{-9} & \text{-}   & \text{-} \\
20 & 8.94\times 10^{-2} & 1.52\times 10^{-5} & 4.56\times 10^{-9} & 5.33\times 10^{-15} & \text{-} \\
25 & 8.88\times 10^{-2} & 1.50\times 10^{-5} & 4.53\times 10^{-9} & 2.78\times 10^{-15} & 5.55\times 10^{-16}\\
30 & 8.85\times 10^{-2} & 1.49\times 10^{-5} & 4.51\times 10^{-9} & 2.73\times 10^{-15} & 5.55\times 10^{-16}\\  \noalign{\hrule height 1pt}
\end{array}\]
\end{table}

Tables \ref{tab:p2_cheb} and \ref{tab:p2_leg} show that the difference in error between Chebyshev and Legendre polynomials for Problem \#2 is small: the Chebyshev polynomials perform slightly better than the Legendre polynomials as the number of points in the domain and number of basis functions increase. 

\section{Conclusions}\label{sec:Conclusions}
This article illustrated that the structure of the univariate TFC \ce\ is composed of a free function and constraint terms, which contain products of projection functionals and switching functions. A method to calculate the projection functionals and switching functions was demonstrated, and their properties were defined. Then, these properties were used as the foundation for mathematical proofs related to the univariate \ce. 

In addition, the projection/switching perspective of the univariate \ce\ lead directly to a multivariate extension via recursive application of the univariate theory. Since these multivariate \ces\ where built from univariate \ces, it was fairly simple to extend the mathematical proofs to the multivariate case as well. In the end, it was concluded that the univariate and multivariate TFC \ces\ are surjective, projection functionals whose domain is the space of all real-valued functions that are defined at the embedded constraints and whose codomain is the space of all real-valued functions that satisfy said constraints. Additionally, a method for compactly writing multivariate \ces\ via tensors was provided. 

After introducing the multivariate TFC in this way, a methodology for solving PDEs via the TFC was presented. This methodology included choosing the free function to be a linear combination of multivariate orthogonal polynomials, discretizing the domain via collocation points, and finally minimizing the residual of the PDE via least-squares. Two example PDEs solved using this methodology were presented, and when available, the TFC solution accuracy was compared with other state-of-the-art methods. 

While this article focused on using orthogonal basis functions, namely, Chebyshev and Legendre orthogonal polynomials, as the free function $g(\B{x})$, and ultimately a linear/nonlinear least-square technique to find the unknown parameters $\B{\xi}$, the technique is not limited to this. At its heart, the TFC approach is a way to derive functionals which analytically satisfy the specified constraints. In other words, when solving ODEs or PDEs, these functionals transform a constrained optimization problem into an unconstrained optimization problem, and therefore a myriad valid definitions of $g(\B{x})$ and optimization schemes exist. For example, deep neural networks, support vector machines, and extreme learning machines have been used as free functions in the past. These other free function choices and associated optimization schemes are useful, as for sufficiently complex problems, the number of multivariate basis functions needed to estimate a PDE with sufficient accuracy can become computationally prohibitive.

As defined in this article, the TFC multivariate \ces\ are capable of embedding value constraints, derivative constraints, and linear combinations thereof. However, other constraint types such as integral, component, and inequality constraints were not discussed. Future work will focus on incorporating these other constraint types. In addition, a more in-depth comparison between TFC and other state-of-the-art methods on a variety of PDEs is likely forthcoming.

\vspace{6pt}
\authorcontributions{Conceptualization, C.L. and H.J.; Formal analysis, C.L. and H.J.; Methodology, C.L. and H.J.; Software, C.L. and H.J.; Supervision, D.M.; Validation, C.L. and H.J.; Writing---original draft, C.L. and H.J.; Writing---review \& editing, C.L., H.J. and D.M.}

\funding{This work was supported by a NASA Space Technology Research Fellowship, Leake [NSTRF 2019] Grant \#: 80NSSC19K1152 and Johnston [NSTRF 2019] Grant \#: 80NSSC19K1149.}


\conflictsofinterest{The authors declare no conflict of interest.} 

\newpage
\abbreviations{The following abbreviations are used in this manuscript:\\

\noindent 
\begin{tabular}{@{}ll}
BVP & Boundary Value Problem\\
FEM & Finite Element Method\\
ODE & Ordinary Differential Equation\\
PDE & Partial Differential Equation\\
TFC & Theory of Functional Connections \\
X-TFC & Extreme Theory of Functional Connections
\end{tabular}}

\appendixtitles{yes} 
\appendix

\section{Orthogonal Polynomials\label{appx:OP}}
\unskip
\subsection{Chebyshev Orthogonal Polynomials}

Chebyshev orthogonal polynomials are two sets of basis functions, the first and the second kind. They are usually indicated as $T_k (z)$ and $U_k (z)$, respectively. This subsection summarizes the main properties of the first kind, $T_k (z)$, only, which are defined on the domain $z\in[-1,+1]$ and with the measure $\dd \mu (z) = \dfrac{1}{\sqrt{1 - z^2}} \dd z$. These polynomials can be generated using the following useful recursive function,
\begin{equation*}
    T_{k + 1} = 2 \, z \, T_k - T_{k - 1} \qquad \text{starting from:} \; \begin{cases} T_0 = 1, \\ T_1 = z.\end{cases}
\end{equation*}

\noindent Also, all the derivatives of Chebyshev orthogonal polynomials can be computed in a recursive way, starting from
\begin{equation*}
    \dfrac{\dd T_0}{\dd z} = 0, \quad \dfrac{\dd T_1}{\dd z} = 1 \qquad \text{and} \qquad \dfrac{\dd^d T_0}{\dd z^d} = \dfrac{\dd^d T_1}{\dd z^d} = 0 \quad (\forall \; d > 1),
\end{equation*}
while the subsequent derivatives of $T_{k + 1} (z)$ can be derived directly from the recursive definition. For~$k \ge 1$ they are,
\begin{equation*}
    \begin{array}{ccccc}
        \dfrac{\dd T_{k+1}}{\dd z} &=& 2 \, \left(T_k + z \, \dfrac{\dd T_k}{\dd z}\right) &- \dfrac{\dd T_{k-1}}{\dd z} \\ [8pt]
        \dfrac{\dd^2 T_{k+1}}{\dd z^2} &=& 2 \left(2 \, \dfrac{\dd T_k}{\dd z} + z \, \dfrac{\dd^2 T_k}{\dd z^2}\right) &- \dfrac{\dd^2 T_{k-1}}{\dd z^2} \\ [4pt]
        \vdots & ~ & \vdots & \vdots \\ [4pt]
        \dfrac{\dd^d T_{k+1}}{\dd z^d} &=& 2 \left( d \, \dfrac{\dd^{d-1} T_k}{\dd z^{d-1}} + z \, \dfrac{\dd^d T_k}{\dd z^d}\right) &- \dfrac{\dd^d T_{k-1}}{\dd z^d}; & (\forall \; d \ge 1).
    \end{array}
\end{equation*}

\noindent The inner product of two Chebyshev orthogonal polynomials satisfies the orthogonality property,
\begin{equation*}
    \langle T_i (z), T_j (z)\rangle = \int_{-1}^{+1} T_i (z) \, T_j (z) \, \dfrac{1}{\sqrt{1 - z^2}} \dd z = \begin{cases} = 0 & {\rm if} \quad i \ne j \\ = \pi & {\rm if} \quad i = j = 0 \\ = \pi/2 & {\rm if} \quad i = j \ne 0\end{cases}.
\end{equation*}
where the orthogonality appears when taking two distinct Chebyshev orthogonal polynomials with indices, $i \ne j$.

\subsection{Legendre Orthogonal Polynomials}

The Legendre orthogonal polynomials, $L_k (z)$, are also defined on the domain $z\in[-1, +1]$, with measure $\dd \mu (z) = \dd z$. These polynomials can also be generated recursively by,
\begin{equation}\label{eq:legendre_polynomials}
    L_{k+1} = \dfrac{2k+1}{k+1} \, z \, L_k - \dfrac{k}{k+1} \, L_{k-1} \qquad \text{starting with:} \; \begin{cases} L_0 = 1 \\ L_1 = z\end{cases}
\end{equation}

\noindent All derivatives of Legendre orthogonal polynomials can be computed in a recursive way, starting from,
\begin{equation*}
    \dfrac{\dd L_0}{\dd z} = 0, \quad \dfrac{\dd L_1}{\dd z} = 1 \qquad \text{and} \qquad \dfrac{\dd^d L_0}{\dd z^d} = \dfrac{\dd^d L_1}{\dd z^d} = 0 \quad (\forall \; d > 1),
\end{equation*}
while the subsequent derivatives of Equation (\ref{eq:legendre_polynomials}) for $k \ge 1$ can be computed in cascade,
\begin{equation*}
    \begin{array}{ccccc}
        \dfrac{\dd L_{k+1}}{\dd z} &=& \dfrac{2k+1}{k+1} \left(L_k + z \dfrac{\dd L_k}{\dd z}\right) & - \dfrac{k}{k+1} \dfrac{\dd L_{k-1}}{\dd z} \\ [8pt]
        \dfrac{\dd^2 L_{k+1}}{\dd z^2} &=& \dfrac{2k+1}{k+1} \left(2\dfrac{\dd L_k}{\dd z} + z \dfrac{\dd^2 L_k}{\dd z^2}\right) & - \dfrac{k}{k+1} \dfrac{\dd^2 L_{k-1}}{\dd z^2} \\ [4pt]
        \vdots & ~ & \vdots & \vdots \\ [4pt]
        \dfrac{\dd^d L_{k+1}}{\dd z^d} &=& \dfrac{2k+1}{k+1} \left(d\dfrac{\dd^{d-1} L_k}{\dd z^{d-1}} + z \dfrac{\dd^d L_k}{\dd z^d}\right) & - \dfrac{k}{k+1} \dfrac{\dd^d L_{k-1}}{\dd z^d}; & (\forall \; d \ge 1).
    \end{array}
\end{equation*}

\noindent Additionally, the orthogonality of Legendre polynomials is given by,
\begin{equation*}
    \langle L_i (z), L_j (z)\rangle = \int_{-1}^{+1} L_i (z) \, L_j (z) \, \dd z = \dfrac{2}{2 i + 1} \, \delta_{ij} \qquad \text{where} \qquad \delta_{ij} = \begin{cases} 1, \quad i = j\\ 0, \quad i\ne j\end{cases}.
\end{equation*}

\section{Solution Times  \label{appx:data}}
\unskip
\subsection{Problem \#1 Solution Times}
\vspace{-6pt}


\begin{table}[H]
\caption{Least-squares solution time using Chebyshev polynomials in milliseconds for Problem \#1.} 
\[\begin{array}{c|ccccc}
\noalign{\hrule height 1.0pt}

\tikz{\node[below left, inner sep=1pt] (n) {\textbf{n}};%
      \node[above right,inner sep=1pt] (m) {\textbf{m}};%
      \draw (n.north west|-m.north west)--(n.south east-|m.south east);}
 & \textbf{5} & \textbf{10} & \textbf{15} & \textbf{20} & \textbf{25}\\
\hline
5  & 0.46 & \text{-} & \text{-} & \text{-} & \text{-} \\
10 & 0.49 & 1.45 & \text{-} & \text{-} & \text{-} \\
15 & 0.61 & 1.88 & 4.90 & \text{-} & \text{-} \\
20 & 0.68 & 2.42 & 6.50 & 13.83 & \text{-} \\
25 & 0.92 & 3.93 & 8.23 & 18.07 & 35.63\\
30 & 0.97 & 4.09 & 13.23 & 22.09 & 50.88\\ 
\noalign{\hrule height 1.0pt}
\end{array}\]
\end{table}
\unskip

\begin{table}[H]
\caption{Total solution time using Chebyshev polynomials in seconds for Problem \#1.}
\[\begin{array}{c|ccccc}
\noalign{\hrule height 1.0pt}

\tikz{\node[below left, inner sep=1pt] (n) {\textbf{n}};%
      \node[above right,inner sep=1pt] (m) {\textbf{m}};%
      \draw (n.north west|-m.north west)--(n.south east-|m.south east);}
 & \textbf{5} & \textbf{10} & \textbf{15} & \textbf{20} & \textbf{25}\\
\hline
5  & 0.053 & \text{-} & \text{-} & \text{-} & \text{-} \\
10 & 0.119 & 0.128 & \text{-} & \text{-} & \text{-} \\
15 & 0.210 & 0.288 & 0.430 & \text{-} & \text{-} \\
20 & 0.399 & 0.584 & 0.824 & 1.114 & \text{-} \\
25 & 0.697 & 1.011 & 1.582 & 2.344 & 3.609\\
30 & 1.045 & 1.871 & 2.909 & 4.683 & 8.164\\  \noalign{\hrule height 1.0pt}
\end{array}\]
\end{table}
\unskip

\begin{table}[H]
\caption{Least-squares solution time using Legendre polynomials in milliseconds for Problem \#1.}
\[\begin{array}{c|ccccc}
\noalign{\hrule height 1.0pt}

\tikz{\node[below left, inner sep=1pt] (n) {\textbf{n}};%
      \node[above right,inner sep=1pt] (m) {\textbf{m}};%
      \draw (n.north west|-m.north west)--(n.south east-|m.south east);}
 & \textbf{5} & \textbf{10} & \textbf{15} & \textbf{20} & \textbf{25}\\
\hline
5  & 0.77 & \text{-} & \text{-} & \text{-} & \text{-} \\
10 & 0.58 & 1.41 & \text{-} & \text{-} & \text{-} \\
15 & 0.56 & 1.82 & 4.70 & \text{-} & \text{-} \\
20 & 0.65 & 2.56 & 6.62 & 13.74 & \text{-} \\
25 & 0.86 & 3.15 & 7.93 & 17.71 & 35.61\\
30 & 1.32 & 3.53 & 9.47 & 22.06 & 46.26\\ \noalign{\hrule height 1.0pt}
\end{array}\]
\end{table}
\unskip

\begin{table}[H]
\caption{Total solution time using Legendre polynomials in seconds for Problem \#1.}
\[\begin{array}{c|ccccc}
\noalign{\hrule height 1.0pt}

\tikz{\node[below left, inner sep=1pt] (n) {\textbf{n}};%
      \node[above right,inner sep=1pt] (m) {\textbf{m}};%
      \draw (n.north west|-m.north west)--(n.south east-|m.south east);}
 & \textbf{5} & \textbf{10} & \textbf{15} & \textbf{20} & \textbf{25}\\
\hline
5  & 0.061 & \text{-} & \text{-} & \text{-} & \text{-} \\
10 & 0.111 & 0.164 & \text{-} & \text{-} & \text{-} \\
15 & 0.214 & 0.298 & 0.368 & \text{-} & \text{-} \\
20 & 0.394 & 0.613 & 0.863 & 1.108 & \text{-} \\
25 & 0.692 & 1.095 & 1.606 & 2.389 & 3.505\\
30 & 1.031 & 1.689 & 2.765 & 4.696 & 7.063\\  \noalign{\hrule height 1.0pt}
\end{array}\]
\end{table}

\subsection{Problem \#2 Solution Times}
\vspace{-6pt}

\begin{table}[H]
\caption{Least-squares solution time using Chebyshev polynomials in milliseconds for Problem \#2.}
\[\begin{array}{c|ccccc}
\noalign{\hrule height 1.0pt}

\tikz{\node[below left, inner sep=1pt] (n) {\textbf{n}};%
      \node[above right,inner sep=1pt] (m) {\textbf{m}};%
      \draw (n.north west|-m.north west)--(n.south east-|m.south east);}
 & \textbf{5} & \textbf{10} & \textbf{15} & \textbf{20} & \textbf{25}\\
\hline
5  & 1.76 & \text{-} & \text{-} & \text{-} & \text{-} \\
10 & 2.14 & 5.56 & \text{-} & \text{-} & \text{-} \\
15 & 2.38 & 7.79 & 27.30 & \text{-} & \text{-} \\
20 & 2.99 & 10.17 & 36.38 & 71.75 & \text{-} \\
25 & 3.90 & 13.00 & 42.40 & 93.11 & 144.5\\
30 & 4.62 & 18.20 & 48.68 & 105.7 & 185.8\\  \noalign{\hrule height 1.0pt}
\end{array}\]
\end{table}
\unskip

\begin{table}[H]
\caption{Total solution time using Chebyshev polynomials in seconds for Problem \#2.}
\[\begin{array}{c|ccccc}
\noalign{\hrule height 1.0pt}

\tikz{\node[below left, inner sep=1pt] (n) {\textbf{n}};%
      \node[above right,inner sep=1pt] (m) {\textbf{m}};%
      \draw (n.north west|-m.north west)--(n.south east-|m.south east);}
 & \textbf{5} & \textbf{10} & \textbf{15} & \textbf{20} & \textbf{25}\\
\hline
5  & 0.547 & \text{-} & \text{-} & \text{-} & \text{-} \\
10 & 1.905 & 1.384 & \text{-} & \text{-} & \text{-} \\
15 & 4.045 & 3.337 & 4.460 & \text{-} & \text{-} \\
20 & 7.460 & 6.753 & 9.771 & 10.97 & \text{-} \\
25 & 11.93 & 12.95 & 17.43 & 21.15 & 22.88\\
30 & 18.94 & 21.98 & 28.90 & 38.92 & 45.82\\  \noalign{\hrule height 1.0pt}
\end{array}\]
\end{table}
\unskip

\begin{table}[H]
\caption{Least-squares solution time using Legendre polynomials in milliseconds for Problem \#2.}
\[\begin{array}{c|ccccc}
\noalign{\hrule height 1.0pt}

\tikz{\node[below left, inner sep=1pt] (n) {\textbf{n}};%
      \node[above right,inner sep=1pt] (m) {\textbf{m}};%
      \draw (n.north west|-m.north west)--(n.south east-|m.south east);}
 & \textbf{5} & \textbf{10} & \textbf{15} & \textbf{20} & \textbf{25}\\
\hline
5  & 1.64 & \text{-} & \text{-} & \text{-} & \text{-} \\
10 & 1.90 & 5.46 & \text{-} & \text{-} & \text{-} \\
15 & 2.24 & 10.63 & 21.52 & \text{-} & \text{-} \\
20 & 3.27 & 11.74 & 29.56 & 69.51 & \text{-} \\
25 & 3.79 & 13.82 & 36.32 & 91.89 & 145.1\\
30 & 4.17 & 13.73 & 51.42 & 112.4 & 181.4\\  \noalign{\hrule height 1.0pt}
\end{array}\]
\end{table}
\unskip

\begin{table}[H]
\caption{Total solution time using Legendre polynomials in seconds for Problem \#2.}
\[\begin{array}{c|ccccc}
\noalign{\hrule height 1.0pt}

\tikz{\node[below left, inner sep=1pt] (n) {\textbf{n}};%
      \node[above right,inner sep=1pt] (m) {\textbf{m}};%
      \draw (n.north west|-m.north west)--(n.south east-|m.south east);}
 & \textbf{5} & \textbf{10} & \textbf{15} & \textbf{20} & \textbf{25}\\
\hline
5  & 0.508 & \text{-} & \text{-} & \text{-} & \text{-} \\
10 & 1.684 & 1.330 & \text{-} & \text{-} & \text{-} \\
15 & 3.875 & 4.489 & 3.849 & \text{-} & \text{-} \\
20 & 7.284 & 7.882 & 8.078 & 10.56 & \text{-} \\
25 & 11.63 & 12.86 & 15.50 & 21.35 & 23.15\\
30 & 17.57 & 18.44 & 28.46 & 40.53 & 45.60\\  \noalign{\hrule height 1.0pt}
\end{array}\]
\end{table}

\reftitle{References}

\begin{thebibliography}{-------}
\providecommand{\natexlab}[1]{#1}

\bibitem[Mortari(2017{\natexlab{a}})]{U-TFC}
Mortari, D.
\newblock The Theory of Connections: Connecting Points.
\newblock {\em Mathematics} {\bf 2017}, {\em 5}.
\newblock
  doi:{\changeurlcolor{black}\href{https://doi.org/10.3390/math5040057}{\detokenize{10.3390/math5040057}}}.

\bibitem[Mortari(2017{\natexlab{b}})]{LDE}
Mortari, D.
\newblock Least-Squares Solution of Linear Differential Equations.
\newblock {\em Mathematics} {\bf 2017}, {\em 5}, 48.
\newblock
  doi:{\changeurlcolor{black}\href{https://doi.org/10.3390/math5040048}{\detokenize{10.3390/math5040048}}}.

\bibitem[Mortari \em{et~al.}(2019)Mortari, Johnston, and Smith]{NDE}
Mortari, D.; Johnston, H.; Smith, L.
\newblock High Accuracy Least-squares Solutions of Nonlinear Differential
  Equations.
\newblock {\em J. Comput. Appl.~Math.} {\bf 2019},
  {\em 352},~293--307.
\newblock
  doi:10.1016/j.cam.2018.12.007.

\bibitem[Johnston and Mortari(2019)]{hybrid_tfc}
Johnston, H.; Mortari, D.
\newblock Least-squares Solutions of Boundary-value Problems in Hybrid Systems.  \emph{arXiv}
   \textbf{2019},  arXiv:math.OC/1911.04390.

\bibitem[Furfaro and Mortari(2019)]{EOL_EOI}
Furfaro, R.; Mortari, D.
\newblock {Least-squares Solution of a Class of Optimal Space Guidance Problems
  via Theory of Connections}.
\newblock {\em ACTA~Astronaut.} {\bf 2019}.
\newblock
  doi:10.1016/j.actaastro.2019.05.050.

\bibitem[Johnston \em{et~al.}(2020)Johnston, Schiassi, Furfaro, and
  Mortari]{FOL}
Johnston, H.; Schiassi, E.; Furfaro, R.; Mortari, D.
\newblock Fuel-Efficient Powered Descent Guidance on Large Planetary Bodies via
  Theory of Functional Connections.  \emph{arXiv}  \textbf{2020},
  arXiv:math.OC/2001.03572.

\bibitem[Mai and Mortari()]{QP_NLP}
Mai, T.; Mortari, D.
\newblock {Theory of Functional Connections Applied to Nonlinear Programming
  under Equality Constraints}.
\newblock In~Proceeding of the 2019 AAS/AIAA Astrodynamics Specialist
  Conference, Portland, ME, USA, 11--15  August  2019; Paper AAS 19-675 . 

\bibitem[Johnston \em{et~al.}(2019)Johnston, Leake, Efendiev, and
  Mortari]{Selected}
Johnston, H.; Leake, C.; Efendiev, Y.; Mortari, D.
\newblock Selected Applications of the Theory of Connections: A~Technique for
  Analytical Constraint Embedding.
\newblock {\em Mathematics} {\bf 2019}, {\em 7}.
\newblock
  doi:{\changeurlcolor{black}\href{https://doi.org/10.3390/math7060537}{\detokenize{10.3390/math7060537}}}.

\bibitem[Mortari and Leake(2019)]{M-TFC}
Mortari, D.; Leake, C.
\newblock The Multivariate Theory of Connections.
\newblock {\em Mathematics} {\bf 2019}, {\em 7}.
\newblock
  doi:{\changeurlcolor{black}\href{https://doi.org/10.3390/math7030296}{\detokenize{10.3390/math7030296}}}.

\bibitem[Leake \em{et~al.}(2019)Leake, Johnston, Smith, and Mortari]{SVM}
Leake, C.; Johnston, H.; Smith, L.; Mortari, D.
\newblock Analytically Embedding Differential Equation Constraints into Least
  Squares Support Vector Machines Using the Theory of Functional Connections.
\newblock {\em Mach. Learn. Knowl. Extr.} {\bf 2019}, {\em
  1},~1058--1083.
\newblock
  doi:{\changeurlcolor{black}\href{https://doi.org/10.3390/make1040060}{\detokenize{10.3390/make1040060}}}.

\bibitem[Leake and Mortari(2020)]{DeepTfc}
Leake, C.; Mortari, D.
\newblock {Deep Theory of Functional Connections: A New Method for Estimating
  the Solutions of Partial Differential Equations}.
\newblock {\em Mach. Learn. Knowl. Extr.} {\bf 2020}, {\em
  2},~37--55.

\bibitem[Schiassi \em{et~al.}(2020)Schiassi, Leake, Florio, Johnston, Furfaro,
  and Mortari]{XTFC}
Schiassi, E.; Leake, C.; Florio, M.D.; Johnston, H.; Furfaro, R.; Mortari, D.
\newblock Extreme Theory of Functional Connections: A Physics-Informed Neural
  Network Method for Solving Parametric Differential Equations.  \emph{arXiv}  \textbf{2020},
  arXiv:cs.LG/2005.10632.

\bibitem[Leake and Mortari(2019)]{M-TFC-PDE}
Leake, C.; Mortari, D.
\newblock {An Explanation and Implementation of Multivariate Theory of
  Functional Connections via Examples}.
\newblock  In~Proceeding of the  {AIAA/AAS} Astrodynamics Specialist Conference,  Portland, ME, USA, 11--15  August  2019 

\bibitem[Ye \em{et~al.}(2014)Ye, Gao, Wang, Cheng, Wang, and Sun]{Ye}
Ye, J.; Gao, Z.; Wang, S.; Cheng, J.; Wang, W.; Sun, W.
\newblock {Comparative Assessment of Orthogonal Polynomials for Wavefront
  Reconstruction over the Square Aperture}.
\newblock {\em J.~Opt. Soc. Am.~A} {\bf 2014}, {\em
  31},~2304--2311.

\bibitem[Dunkl and Xu(2014)]{Dunkl}
Dunkl, C.F.; Xu, Y.
\newblock {\em Orthogonal Polynomials of Several Variables}, 2 ed.;
  Encyclopedia of Mathematics and its Applications;  Cambridge University Press: Cambridge, UK, 
   2014.
\newblock
  doi:{\changeurlcolor{black}\href{https://doi.org/10.1017/CBO9781107786134}{\detokenize{10.1017/CBO9781107786134}}}.

\bibitem[Xu(1994)]{Xu}
Xu, Y.
\newblock {Multivariate Orthogonal Polynomials and Operator Theory}.
\newblock {\em Trans. Am. Math. Soc.} {\bf 1994},
  {\em 343},~193--202.

\bibitem[Langtangen(2003)]{nDbasisFunctions}
Langtangen, H.P.
\newblock {\em Computational {Partial} {Differential} {Equations}: {Numerical}
  {Methods} and {Diffpack} {Programming}}; Springer: Berlin/Heidelberg, Germany, 2003; 
\newblock OCLC: 851766084.

\bibitem[Lanczos(1957)]{Collo1}
Lanczos, C., Applied Analysis.
\newblock In {\em Progress in Industrial Mathematics at ECMI 2008}; Dover
  Publications, Inc.: New York, NY, USA,  1957; p. 504.

\bibitem[Wright(1964)]{Collo2}
Wright, K.
\newblock {Chebyshev Collocation Methods for Ordinary Differential Equations.}
\newblock {\em   Comput. J.} {\bf 1964}, {\em 6},~358--365.


\bibitem[Maclaurin \em{et~al.}(2013)Maclaurin, Duvenaud, Johnson, and
  Townsend]{autograd}
Maclaurin, D.; Duvenaud, D.; Johnson, M.; Townsend, J.
\newblock Autograd. 2013. Available online: 
\newblock \url{https://github.com/HIPS/autograd} (accessed on 1 July 2020). 

\bibitem[Baydin \em{et~al.}(2018)Baydin, Pearlmutter, Radul, and
  Siskind]{autodiff}
Baydin, A.G.; Pearlmutter, B.A.; Radul, A.A.; Siskind, J.M.
\newblock Automatic Differentiation in Machine Learning: A Survey.
\newblock {\em J.~Mach. Learn.~Res.} {\bf 2018}, {\em
  18},~1--43.

\bibitem[Lagaris \em{et~al.}(1998)Lagaris, Likas, and Fotiadis]{OrigOdePde}
Lagaris, I.E.; Likas, A.; Fotiadis, D.I.
\newblock {Artificial neural networks for solving ordinary and partial
  differential equations}.
\newblock {\em IEEE Trans. Neural Netw.} {\bf 1998}, {\em
  9},~987--1000.
\newblock
  doi:{\changeurlcolor{black}\href{https://doi.org/10.1109/72.712178}{\detokenize{10.1109/72.712178}}}.

\bibitem[Mall and Chakraverty(2017)]{CNN}
Mall, S.; Chakraverty, S.
\newblock {Single Layer Chebyshev Neural Network Model for Solving Elliptic
  Partial Differential Equations}.
\newblock {\em Neural Process. Lett.} {\bf 2017}, {\em 45},~825--840.
\newblock
  doi:{\changeurlcolor{black}\href{https://doi.org/10.1007/s11063-016-9551-9}{\detokenize{10.1007/s11063-016-9551-9}}}.

\bibitem[Sun \em{et~al.}(2019)Sun, Hou, Yang, Zhang, Weng, and Han]{BNN}
Sun, H.; Hou, M.; Yang, Y.; Zhang, T.; Weng, F.; Han, F.
\newblock {Solving Partial Differential Equation Based on Bernstein Neural
  Network and Extreme Learning Machine Algorithm}.
\newblock {\em Neural Process. Lett.} {\bf 2019}, {\em 50},~1153--1172.
\newblock
  doi:{\changeurlcolor{black}\href{https://doi.org/10.1007/s11063-018-9911-8}{\detokenize{10.1007/s11063-018-9911-8}}}.

\bibitem[Huang \em{et~al.}(2006)Huang, Zhu, and Siew]{ELM}
Huang, G.B.; Zhu, Q.Y.; Siew, C.K.
\newblock { Extreme learning machine: Theory and applications}.
\newblock {\em Neurocomputing} {\bf 2006}, {\em 70},~489--501.
\newblock
  doi:{\changeurlcolor{black}\href{https://doi.org/10.1016/j.neucom.2005.12.126}{\detokenize{10.1016/j.neucom.2005.12.126}}}.

\end{thebibliography}

\end{document}